\documentclass{article}
\usepackage{graphicx} 
\usepackage{amssymb}
\usepackage{amsthm}
\usepackage{aliascnt}

\usepackage[numbers,sort&compress]{natbib}
\bibpunct[, ]{[}{]}{,}{n}{,}{,}
\renewcommand\bibfont{\fontsize{10}{12}\selectfont}

\usepackage{amsmath} 
\usepackage{caption}
\usepackage{algorithmic}
\usepackage[ruled,linesnumbered,vlined]{algorithm2e}
\usepackage{hyperref}
\usepackage{multirow}
\usepackage{subcaption}
\newdimen\figrasterwd
\figrasterwd\textwidth
\usepackage[shortlabels]{enumitem}
\usepackage{xcolor}
\usepackage{cleveref}
\usepackage{enumitem}
\usepackage{geometry}
\usepackage{tikz}
\usetikzlibrary{tikzmark,calc,decorations.pathreplacing}
\usepackage{tikzpagenodes}
\usepackage{float}

\theoremstyle{plain}
\newtheorem{theorem}{Theorem}[section]

\newaliascnt{lemma}{theorem}
\newtheorem{lemma}[lemma]{Lemma}
\aliascntresetthe{lemma}

\newaliascnt{corollary}{theorem}

\aliascntresetthe{corollary}

\newaliascnt{proposition}{theorem}
\newtheorem{proposition}[proposition]{Proposition}
\aliascntresetthe{proposition}

\theoremstyle{definition}
\newaliascnt{definition}{theorem}
\newtheorem{definition}[definition]{Definition}
\aliascntresetthe{definition}

\newaliascnt{example}{theorem}
\newtheorem{example}[example]{Example}
\aliascntresetthe{example}

\newtheorem{assumption}{Assumption}

\theoremstyle{remark}
\newaliascnt{remark}{theorem}
\newtheorem{remark}[remark]{Remark}
\aliascntresetthe{remark}

\newcommand{\K}{\mathcal{K}}

\usepackage{booktabs}      
\usepackage{arydshln}
\DeclareMathOperator{\diag}{diag}

\usepackage{tikz}
\usetikzlibrary{calc}
\usepackage{pgfplots}
\pgfplotsset{compat=newest}
\usetikzlibrary{arrows,arrows.meta,intersections,shapes.geometric,positioning,calc,patterns}

\title{A constraint dissolving method for generally constrained optimization over a nonconvex set\thanks{The research of Pengjie Liu was supported by the National Natural Science Foundation of China (Grant No. 12601615), the China Postdoctoral Science Foundation (Grant No. 2026M793421), and the Jiangsu Funding Program for Excellent Postdoctoral Talent. The research of Xiaoxi Jia was supported by the European Social Fund Plus (ESF Plus) and the Free State of Saxony through the project MORE-KIBA (Project No. 100779486).}}
\author{
Pengjie Liu$^{\dagger}$,
Haonan Shen$^{\dagger}$,
Xiaoxi Jia$^{\ddagger}$\\[1em]
$^{\dagger}$School of Mathematics, China University of Mining and Technology, Xuzhou 221116, China\\[0.3em]
$^{\ddagger}$Faculty of Mathematics, Chemnitz University of Technology, Germany
}

\date{3 Nov. 2026}
\begin{document}

\maketitle

\begin{abstract}
This manuscript studies a constraint dissolving method for solving optimization problems with general constraints, comprising an abstract constraint $G(x)\in C$ and a potentially nonconvex geometric constraint $x\in D$, which, to the best of our knowledge, is the first constraint dissolving type method applicable to this general setting. More specifically, we introduce an auxiliary variable $s$ to decouple $G(x)\in C$ into $G(x)=s$ and $s\in C$, and construct the constraint dissolving mapping acting jointly on both variables. We then propose a penalty-based constraint dissolving method, in which the penalty parameter is updated adaptively according to the feasibility error. Under mild conditions, every accumulation point is Mordukhovich-stationary. By a specific choice of the constraint dissolving mapping, every accumulation point is guaranteed to be feasible, which stands in contrast to standard penalty-type methods in nonconvex settings. Numerical experiments on Rosenbrock, portfolio, and MAXCUT problems demonstrate the effectiveness of the proposed method against the augmented Lagrangian and the penalty decomposition methods.
\end{abstract}


\section{Introduction}
In this paper, we consider the following optimization problem with a geometric constraint:
\begin{equation}\label{Eq:P}\tag{P}
\begin{aligned}
    \min_{x}\quad & f(x)\\
    \text{s.t.}\quad & G(x)\in C, \  x\in D,
\end{aligned}
\end{equation}
where $f:\mathbb R^n\to\mathbb R$ and $G:\mathbb R^n\to\mathbb R^m$ are continuously differentiable. The set $C\subset\mathbb R^m$ is nonempty, closed, and convex, whereas $D\subset\mathbb R^n$ is nonempty and closed, but not necessarily convex. We refer to $x\in D$ as a geometric constraint, allowing $D$ to represent potentially nonconvex structures, such as cardinality, complementarity, either-or, and low-rank constraints.

To the best of our knowledge, penalty-based methods are potentially efficient for solving such types of problems. Specifically, Jia et al. \cite{jia2023augmented} develop an augmented Lagrangian method that handles the geometric constraint in a unified framework, while Kanzow and Lapucci \cite{kanzow2023inexact} propose an inexact penalty decomposition method where an auxiliary variable is employed to separate the geometric constraint $x\in D$. More recently, the constraint dissolving technique, originally introduced for optimization over the Stiefel manifold \cite{xiao2024stiefelpenalty}, has shown great potential for solving \eqref{Eq:P}. This type of method is first developed for problems with equality constraints $G(x)=0$ in \cite{xiao2024dissolving}, and a practical scheme based on Monte Carlo sampling for determining the penalty parameter is provided in \cite{xiao2025cdopt}. It is subsequently generalized to handle equality constraints $G(x)=0$ together with a convex set constraint $x\in D$ in \cite{xiao2025exact}. Beyond these settings, Hu et al. \cite{hu2023improved} transform nonconvex-strongly-convex bilevel optimization problems into unconstrained constraint dissolving problems and develop a unified framework for subgradient-based methods, and the technique has also been applied to feasibility problems over nonconvex sets \cite{xiao2025quadratically}. Building on the works mentioned above, Jia et al. \cite{jia2026constraint} consider this type of method for solving \eqref{Eq:P}, i.e., with the more general constraints $G(x)\in C$ and $x\in D$, where $D$ is possibly nonconvex, and establish theoretical foundations for the constraint dissolving mapping and the corresponding algorithm. However, due to the nonsmoothness of the projection onto $C$, the numerical realization in \cite{jia2026constraint} is restricted to the case $C:=\{0\}$, leaving the general convex case, including standard inequality constraints,  unaddressed. This manuscript aims to resolve this gap, as a substantial extension of \cite{jia2026constraint}. 

To this end, we introduce a slack variable $s\in\mathbb R^m$ to separate the nonlinear coupling from the set constraint $G(x)\in C$ and rewrite
\eqref{Eq:P} as
\begin{equation}\label{Eq:Ps}\tag{$P_s$}
 \begin{aligned}
  \min_{x,s}\quad & f(x)\\
  \text{s.t.}\quad & G(x)-s=0, \ x\in D, \  s\in C.
 \end{aligned}
\end{equation}
This standard reformulation preserves feasible points and local
and global minimizers same as those of \eqref{Eq:P} in the sense that the map
$x\mapsto(x,G(x))$ is bijection; see, e.g.,
\cite[Section~3]{de2023constrained}. Based on this reformulation, the main contributions of this manuscript are threefold.
\begin{itemize}
    \item Theoretically, we introduce a constraint dissolving mapping that acts jointly on the original and slack variables. Extending the analysis of \cite{jia2026constraint} to this reformulation, we prove under a mild constraint qualification that, for every positive penalty parameter, a feasible point $x$ is Mordukhovich stationary for \eqref{Eq:P} if and only if $(x,G(x))$ is Mordukhovich stationary for the resulting constraint dissolving problem. On this basis, we show that every feasible accumulation point of the proposed method is a Mordukhovich stationary point of \eqref{Eq:P}. Moreover, for a specific constraint dissolving mapping, we prove that every accumulation point is feasible for \eqref{Eq:P}, which stands in contrast to classical penalty-type methods in the nonconvex setting, whose accumulation points are not necessarily feasible in general.
    \item Algorithmically, we propose a penalty-based constraint dissolving method, which computes, by a projected gradient method, an approximate Mordukhovich stationary point of the penalty function obtained from constraint dissolving over $(x,s)\in D\times C$. For the constraint dissolving problem to be solved, the constraints $x\in D$ and $s\in C$ are retained, whereas $G(x)=s$ is dissolved into the objective function. Afterwards, the penalty parameter $\beta$ is kept if the feasibility error $\|G(x)-s\|$ has decreased by a prescribed factor and is increased otherwise. This explicit adaptive rule stands in contrast to the fixed penalty parameter in \cite{xiao2025exact} and the nondecreasing update without an explicit rule in \cite{jia2026constraint}.
    \item Numerically, we implement the proposed method with the explicit constraint dissolving mapping constructed in this manuscript and test it on complementarity-constrained Rosenbrock, sparse portfolio, and MAXCUT problems. In accordance with the general setting considered here, these problems include instances that combine a nonconvex set $D$ with a constraint set $C$ that is not a singleton. In the presented experiments, it requires, in most cases, fewer outer iterations than the safeguarded augmented Lagrangian method \cite{jia2023augmented} and the inexact penalty decomposition method \cite{kanzow2023inexact}, while attaining solutions of comparable or better quality.
\end{itemize}

The rest of the paper is organized as follows. In \Cref{sec:pre}, we recall necessary definitions and properties of variational analysis, establish properties of the constraint dissolving mapping, and characterize the relationships between stationary points of \eqref{Eq:P} and those of the reformulated problem via the constraint dissolving mapping. In \Cref{Sec:alg}, we propose a penalty-based constraint dissolving method for solving \eqref{Eq:P} and establish its convergence analysis. In \Cref{sec:realization}, we construct a specific constraint dissolving mapping and prove that every accumulation point generated by the proposed algorithm is feasible for \eqref{Eq:P}. In \Cref{Sec:Num}, we conduct numerical experiments on several benchmark problems, including Rosenbrock, portfolio optimization, and MAXCUT problems, to demonstrate the effectiveness of the proposed algorithm. Finally, \Cref{Sec:conclusion} concludes the manuscript.

\section{Preliminaries}
\label{sec:pre}

In this section, we first recall the tools from variational analysis \cite{rockafellar1998variational} used throughout the paper in \Cref{Sub:va}.
\Cref{Subs:properties_cd} derives key properties of the constraint dissolving mapping. \Cref{Subs:relation_stationarities} investigates  relationships between stationary points of \eqref{Eq:P} and those of the corresponding penalty problem based on the constraint dissolving mapping.

\subsection{Materials for variational analysis}\label{Sub:va}

For any nonempty, closed, but not necessarily convex set
\(D \subseteq \mathbb{R}^n\), we denote
$P_D(x) := \mathop{\arg\min}_{z \in D} \|z-x\|$ and $\operatorname{dist}_D(x) := \inf_{z \in D} \|z-x\|$ 
as the corresponding projection mapping and the distance function,
respectively. When \(D\) is nonconvex, \(P_D(x)\) may fail to
be a singleton, i.e., the projection onto \(D\) does not need to be unique.

Given \(z\in D\), the tangent cone to \(D\) at \(z\) is defined by
\[
\mathcal{T}_D(z)
:=
\left\{
w\in\mathbb{R}^n
\;\middle|\;
\frac{z^k-z}{\tau^k}\to w
\ 
\text{for some }
z^k \underset{D}{\longrightarrow} z,
\ \tau^k\downarrow 0
\right\},
\]
the normal cone to \(D\) at \(z\) is defined as
\[
\mathcal{N}_D(z)
:=
\left\{
v\in\mathbb{R}^n
\;\middle|\;
\langle v,x-z\rangle
\leq
o(\|x-z\|)
\text{ for }x\in D
\right\},
\]
the limiting normal cone is defined as
\[
\mathcal{N}_D^{\lim}(z)
:=
\left\{
v\in\mathbb{R}^n
\;\middle|\;
\text{there exist sequences }
z^k \underset{D}{\rightarrow} z
\text{ and }
v^k\to v
\text{ with }
v^k\in\mathcal{N}_D(z^k)
\right\}.
\]
We set  $\mathcal{N}_D(z)=\emptyset$ and $\mathcal{N}_D^{\lim}(z)=\emptyset$ if $ z\notin D $. Note that the limiting normal cone is stable in the sense that
\begin{equation}
    \limsup_{z\to\bar z}
    \mathcal{N}_D^{\lim}(z)
    =
    \mathcal{N}_D^{\lim}(\bar z),
    \ \forall \ \bar z\in D.
\end{equation}
When \(D\) is convex, the limiting normal cone coincides with the
corresponding standard normal cone, i.e., for \(z\in D\),
\[
\mathcal{N}_D^{\lim}(z)
=
\mathcal{N}_D(z)
:=
\left\{
v\in\mathbb{R}^n
\;\middle|\;
\langle v,x-z\rangle\leq 0,
\ 
\forall \ x\in D
\right\}.
\]

Let us first specify the stationarity concepts for the original problem.
These provide the reference conditions for the reformulations and the
convergence analysis. In this manuscript, we denote the feasible set of \eqref{Eq:P} by $\mathcal F:=\{x\in D:G(x)\in C\}$. 

\begin{definition}\label{def:m-stationarity}
A point $\bar x\in\mathcal F$ is called an M (Mordukhovich)-stationary point for
\eqref{Eq:P} if there exists $\lambda\in\mathbb R^m$ such that
\begin{equation}\label{eq:m-original}
 0\in\nabla f(\bar x)+\nabla G(\bar x)^\top\lambda
                  +\mathcal N_D^{\lim}(\bar x),
 \  \lambda\in\mathcal N_C(G(\bar x)).
\end{equation}
\end{definition}

\begin{definition}\label{def:am-stationarity}
A point \(\bar{x}\in\mathcal{F}\) is called an 
AM (asymptotically Mordukhovich)-stationary point of
\eqref{Eq:P} if there exist sequences $\{x^k\},\{\varepsilon^k\}\subset\mathbb{R}^n$ and $\{\lambda^k\},\{z^k\}\subset\mathbb{R}^m$
such that 
$x^k\to\bar{x}, \varepsilon^k\to 0, z^k\to 0$,
and for all $k\in \mathbb N$, 
\begin{equation}\label{eq:am-original} 
\varepsilon^k
 \in
\nabla f(x^k)
+
\nabla G(x^k)^\top\lambda^k
+
\mathcal{N}_D^{\lim}(x^k), \ \ 
\lambda^k
 \in
\mathcal{N}_C\bigl(G(x^k)-z^k\bigr).  
\end{equation} 
\end{definition}

Note that M-stationarity obviously implies AM-stationarity. The converse is ensured by the following constraint qualification; see \cite[Section~2.2]{jia2023augmented}.

\begin{definition}\label{def:am-regularity}
A point $\bar x\in\mathcal F$ of \eqref{Eq:P} is AM-regular if
\begin{equation}\label{eq:am-regularity}
 \limsup_{\substack{x\to\bar x\\z\to0}}\mathcal M(x,z)
       \subseteq\mathcal M(\bar x,0),
\end{equation}
where $\mathcal M:\mathbb R^n\times\mathbb R^m\rightrightarrows
\mathbb R^n$ is the set-valued mapping defined by
\[
 \mathcal M(x,z):=\nabla G(x)^T\mathcal N_C(G(x)-z)
                         +\mathcal N_D^{\lim}(x).
\]
\end{definition}

As illustrated in \cite{jia2026constraint,xiao2024stiefelpenalty,xiao2024dissolving,xiao2025exact}, a constraint qualification is typically required for the analysis of constraint dissolving methods. As shown in \cite[Lemma~2.15]{jia2026constraint}, AM-regularity is weaker than the constraint qualification employed in \cite{xiao2024stiefelpenalty,xiao2024dissolving,xiao2025exact}. We now state some assumptions required for our analysis.
\begin{assumption}[Properties of the optimization problem \eqref{Eq:P}]
\label{ass:P}
\leavevmode\par\nobreak
\begin{enumerate}[label=(\roman*), ref=\theassumption(\roman*)]
    \item\label{ass:P-i}
    $f:\mathbb{R}^n\to\mathbb{R}$ and
    $G:\mathbb{R}^n\to\mathbb{R}^m$ are continuously differentiable.

    \item\label{ass:P-ii}
    $C\subseteq\mathbb{R}^m$ is nonempty, closed, and convex.
    $D\subseteq\mathbb{R}^n$ is nonempty, closed, but not necessarily convex.

    \item\label{ass:P-iii}
    AM-regularity, as defined in \Cref{def:am-regularity}, holds at every feasible point.
\end{enumerate}
\end{assumption}

\subsection{Properties of constraint dissolving mappings}\label{Subs:properties_cd}

This subsection focuses on the definition and properties of the constraint dissolving mapping tailored for \eqref{Eq:Ps}. For convenience, we denote $H(x,s):=G(x)-s$, $K:=D\times C$, and $\mathcal F_s:=\{(x,s)\in K\,|\, G(x)=s\}$.

\begin{definition}\label{def:cdm}
A mapping
$
 \mathcal A:
          K\to\mathbb R^{n+m}
$
is called a constraint dissolving mapping for $H(x,s)=0$,
if the following conditions hold.
\begin{enumerate}[label=(\alph*)]
 \item \label{item:Aa}$\mathcal A$ is continuously differentiable.
 \item $\mathcal A(x,s)=(x,s)$ holds for all
       $(x,s)\in\mathcal F_s$.
 \item \label{item:Ac}For all $(x,s)\in\mathcal F_s$,
       $\nabla\mathcal A(x,s)^\top \nabla H(x,s)^\top = 0$. 
\end{enumerate}
Moreover, there exists a continuous mapping
$\mathcal R_{\mathcal A}:K\to\mathbb R^{(n+m)\times(n+m)}$
such that:
\begin{enumerate}[label=(\alph*),start=4]
 \item \label{item:Ad} $\mathcal R_{\mathcal A}(x,s)=0$ for all
       $(x,s)\in\mathcal F_s$;
 \item \label{item:Ae}  for any $(x,s)\in\mathcal F_s$, there exists
       $\omega_{x,s}>0$ such that
       \begin{equation}\label{eq:normal-preservation}
       \bigl(\nabla\mathcal A(y,t)-I_{n+m}
                       -\mathcal R_{\mathcal A}(y,t)\bigr)^\top
                \mathcal N_K^{\lim}(y,t)=\{0\}
       \end{equation}
       for all $(y,t)\in K$ with
       $\|(y,t)-(x,s)\|\leq\omega_{x,s}$.
\end{enumerate}
\end{definition}

\Cref{def:cdm} generalizes the constraint dissolving mappings of
\cite{xiao2024dissolving,xiao2025exact}, which were developed for an
equality constraint $G(x)=0$. In our setting, by contrast, the set $C$ can be an
arbitrary nonempty, closed, and convex set. Moreover, whereas
\cite{xiao2025exact} additionally admits a constraint $x\in \mathcal{X}$ with a
convex set $\mathcal{X}$, we incorporate the geometric constraint $x\in D$, where $D$ is allowed to be nonconvex. Therefore, the corresponding limiting normal cone is employed in \eqref{eq:normal-preservation}.

\begin{proposition}\label{prop:cdm-tangent}
Let $\mathcal A$ be any constraint dissolving mapping satisfying \Cref{def:cdm}.
For every $(x,s)\in\mathcal F_s$ and
$(d_x,d_s)\in\mathcal T_{\mathcal F_s}(x,s)$, one has
$\nabla\mathcal A(x,s) ( d_x, d_s) = ( d_x, d_s)$. 
\end{proposition}
\begin{proof}
For any $(d_x, d_s)\in \mathcal T_{\mathcal F_s}(x,s)$ with $(x,s)\in \mathcal F_s$, there exists a sequence $(x^j,s^j)\in\mathcal F_s$ and $t_j\downarrow0$ such that
$(x^j,s^j)\to(x,s)$ and
\begin{equation}\label{Eq:tangentforA}
\frac{(x^j, s^j)-(x, s)}{t_j}\to(d_x,d_s).
\end{equation}
By differentiability of $\mathcal A$ as in \Cref{def:cdm} \ref{item:Aa}, one has
\begin{equation}\label{Eq:difA}
    \mathcal A(x^j,s^j)-\mathcal A(x,s) = \nabla \mathcal A(x,s)\left((x^j, s^j)-(x,s)\right) + o \left(\|(x^j, s^j)-(x,s)\|\right). 
\end{equation} 
\Cref{def:cdm}(b) gives
$\mathcal A(x^j,s^j)-\mathcal A(x,s)=(x^j, s^j)-(x, s)$.  
Dividing it by $t_j$, then \eqref{Eq:tangentforA} and \eqref{Eq:difA} yield the assertion.
\end{proof}

The following result shows that if $\bar x$ satisfies AM-regularity, then the analogous regularity condition, formulated in terms of the set-valued mapping $(x,s) \mapsto \operatorname{Range}(\nabla H(x,s)^{\top}) + \mathcal{N}_{K}^{\lim}(x,s)$, holds at the pair $(\bar{x}, G(\bar{x}))$.

\begin{lemma}
\label{lem:multiplier-closedness}
If $\bar x\in\mathcal F$ is AM-regular for \eqref{Eq:P}, then one has
$\limsup\limits_{\substack{(x,s)\to(\bar x, G(\bar x))}}
       \hat{\mathcal M} (x,s)
   \subseteq \hat{\mathcal M} (\bar x, G(\bar x))$, 
where 
\begin{equation}\label{Eq:hatM}
    \hat{\mathcal M} (x,s):= \operatorname{Range}(\nabla H(x,s)^{\top}) + \mathcal{N}_{K}^{\lim}(x,s).
\end{equation}
\end{lemma}
\begin{proof}
Let us take any $(a,b)\in \limsup_{\substack{(x,s)\to(\bar x, G(\bar x))}}
       \hat{\mathcal M} (x,s)$, then there exist sequences $\{(x^k, s^k)\}\subset K$ and $\{(a^k, b^k)\}$ satisfying $(x^k, s^k)\to (\bar x, G(\bar x))$ and $(a^k, b^k) \to (a,b)$, such that 
$(a^k,b^k)\in  \hat{\mathcal M} (x^k,s^k) $. For any $k\in \mathbb N$, by \cite[Proposition~6.41]{rockafellar1998variational}, there exist $\lambda^k\in \mathbb R^m$, 
$u^k\in\mathcal N_D^{\lim}(x^k)$ and
$v^k\in\mathcal N_C(s^k)$ such that
$a^k=\nabla G(x^k)^\top \lambda^k+u^k$ and $b^k=-\lambda^k+v^k$. 
Consequently,
\[
 a^k+\nabla G(x^k)^\top b^k
       =\nabla G(x^k)^\top v^k+u^k
       \in\mathcal M(x^k,z^k) \  \text{with} \ z^k:=G(x^k)- s^k.
\]
By the continuity of $G$, we have $z^k \to 0$. Then the continuity of $\nabla G$ and AM-regularity at $\bar x$ yield
$a+\nabla G(\bar x)^\top b\in\mathcal M(\bar x,0)= \nabla G(\bar x)^\top \mathcal N_C (G(\bar x))+ \mathcal N_D^{\text{lim}}(\bar x).$ 
Thus, there exists 
$v\in\mathcal N_C(G(\bar x))$ and
$u\in\mathcal N_D^{\lim}(\bar x)$ such that
$a+ \nabla G(\bar x)^\top b = \nabla G(\bar x)^\top v + u$. 
Let us set $\xi:=v-b\in \mathbb R^m$, then we have
$a = \nabla G(\bar x)^ \top \xi + u$ and $b=-\xi+v$, 
it follows from \cite[Proposition~6.41]{rockafellar1998variational} that $(a, b)\in \hat{\mathcal M}(\bar x, G(\bar x))$. Therefore, the assertion holds.
\end{proof}

\begin{remark}
To be consistent with the form of \Cref{def:am-regularity}, the set-valued mapping associated with the constraints $H(x,s)=0$ and $(x,s)\in K$ would read
\[
  \widetilde{\mathcal M}\big((x,s),z\big)
  := \nabla H(x,s)^{\top}\mathcal N_{\{0\}}\big(H(x,s)-z\big)+\mathcal N_K^{\lim}(x,s),
\]
where the role of $z$ is to ensure that $H(x,s)-z=0$ and consequently $\mathcal N_{\{0\}}\big(H(x,s)-z\big)=\mathbb R^m$. In this way, the auxiliary variable $z$ can be eliminated, thus, for convenience, we adopt the simplified form \eqref{Eq:hatM}.
\end{remark}

By \cite[Theorem~3.16]{mehlitz2020asymptotic}, we directly have the following result.
\begin{lemma}\label{lem:normal-inclusion}
If $ x\in\mathcal F$ is AM-regular, then one has
\begin{equation}\label{eq:normal-inclusion}
 \mathcal N_{\mathcal F_s}^{\lim}( x, G(x))
 \subseteq\operatorname{Range}(\nabla H(x, G(x))^\top)
                    +\mathcal N_K^{\lim}(x, G( x)).
\end{equation}
\end{lemma}

Meanwhile, we have the following results.

\begin{proposition}\label{prop:properties-cdm}
Let Assumption~\ref{ass:P} hold and $\mathcal A$ satisfy
\Cref{def:cdm}. For any $(x,s)\in\mathcal F_s$, one has
\begin{enumerate}[label=(\roman*)]
  \item\label{item:cdm-range}
  $\nabla\mathcal A(x,s)^\top \operatorname{Range}(\nabla H(x,s)^\top) = \{0\}$,
  \item\label{item:cdm-normal}
  $(\nabla\mathcal A(x,s)-I_{n+m})^\top \mathcal N_K^{\lim}(x,s) = \{0\}$,
  \item\label{item:cdm-sum}
  $\nabla\mathcal A(x,s)^\top \bigl(\operatorname{Range}(\nabla H(x,s)^\top) + \mathcal N_K^{\lim}(x,s)\bigr) = \mathcal N_K^{\lim}(x,s)$.
\end{enumerate}
Moreover, for any $w\in\mathbb R^{n+m}$, one has
\begin{equation*}
0\in\nabla\mathcal A(x,s)^\top w+\mathcal N_K^{\lim}(x,s)
\ \Longrightarrow \ 
\nabla\mathcal A(x,s)^\top w-w
    \in\operatorname{Range}(\nabla H(x,s)^\top).
\end{equation*}
\end{proposition}

\begin{proof}
The identities \ref{item:cdm-range} and \ref{item:cdm-normal} follow from \Cref{def:cdm} \ref{item:Ac}-\ref{item:Ae}. The identity \ref{item:cdm-sum} is the direct consequence of \ref{item:cdm-range} and \ref{item:cdm-normal}. 
For any $w\in \mathbb R^{n+m}$ satisfying $0\in\nabla\mathcal A(x,s)^\top w+\mathcal N_K^{\lim}(x,s)$, we have
$-\nabla\mathcal A(x,s)^\top w\in\mathcal N_K^{\lim}(x,s)$. Then 
\ref{item:cdm-normal} yields
\begin{equation}\label{eq:idempotence-on-w}
\bigl(\nabla\mathcal A(x,s)^\top \bigr)^2w
    =\nabla\mathcal A(x,s)^\top w.
\end{equation}
For every $d\in\mathcal T_{\mathcal F_s}(x,s)$,
Proposition~\ref{prop:cdm-tangent} gives
$\bigl\langle d,\nabla\mathcal A(x,s)^\top w-w\bigr\rangle
=
\bigl\langle
(\nabla\mathcal A(x,s)-I_{n+m})d,w
\bigr\rangle
=0$. 
It implies that
\[
\nabla\mathcal A(x,s)^\top w-w
\in \left(\mathcal T_{\mathcal F_s}(x,s) \right)^{\circ}
\subseteq\mathcal N_{\mathcal F_s}^{\lim}(x,s)
\subseteq
\operatorname{Range}(\nabla H(x,s)^\top)
+\mathcal N_K^{\lim}(x,s).
\]
where the last inclusion is from \Cref{lem:normal-inclusion}. Then there exist some $\lambda \in \mathbb R^m$ and $q\in \mathcal N_K^{\text{lim}}(x,s)$ such that
$\nabla\mathcal A(x,s)^\top w-w = \nabla H(x,s)^\top \lambda +q$.  
By \eqref{eq:idempotence-on-w}, \ref{item:cdm-range}, and \ref{item:cdm-normal}, we have $q=0$ and thus 
$\nabla\mathcal A(x,s)^\top w-w = \nabla H(x,s)^\top \lambda \in \operatorname{Range}(\nabla H(x,s)^\top)$. 
\end{proof}

Based on \Cref{prop:properties-cdm}, taking the same argument as the proof in \cite[Theorem~2.8]{jia2026constraint} directly implies the following result. 
\begin{theorem}\label{thm:stationarity-equivalence}
Let Assumption~\ref{ass:P} hold and $\mathcal A$ satisfy
\Cref{def:cdm}. For any
$(x,s)\in\mathcal F_s$ and $w\in\mathbb R^{n+m}$, one has
\begin{equation*} 
 0 \in\nabla\mathcal A(x,s)^\top w+\mathcal N_K^{\lim}(x,s) \ 
 \Longleftrightarrow \ 
 0 \in w+\operatorname{Range}(\nabla H(x,s)^\top)+\mathcal N_K^{\lim}(x,s). 
\end{equation*}
\end{theorem}

\subsection{Relationships between stationarities}\label{Subs:relation_stationarities}

As illustrated above, the constraint dissolving mapping employed in this manuscript is jointly with the original variable $x$ and slack variable $s$. For the convenience of analysis, we partition the constraint dissolving mapping $\mathcal A(x,s)$ with respect to the blocks $x$ and $s$, namely,
\begin{equation*}
    \mathcal A(x, s) = (\mathcal A_x(x,s),\, \mathcal A_s(x,s)),
\end{equation*}
where $\mathcal A_x(x,s)\in \mathbb R^n$ and $\mathcal A_s(x,s)\in \mathbb R^m$. Then, by \Cref{def:cdm}, $\mathcal A_x$ and $\mathcal A_s$ satisfy, in
particular, the following properties:
\begin{enumerate}[label=(\roman*)]
    \item $\mathcal A_x$ and $\mathcal A_s$ are continuously differentiable.
    \item $\mathcal A_x(x,s)=x$ and $\mathcal A_s(x,s)=s$ for all $(x,s)\in \mathcal F_s$.
\end{enumerate}
Moreover, the Jacobian of $\mathcal A(x,s)$ admits
\begin{equation*}
    \nabla\mathcal A(x,s) =
   \begin{bmatrix}
        \nabla_x\mathcal A_x(x,s) & \nabla_s\mathcal A_x(x,s)\\[2pt]
        \nabla_x\mathcal A_s(x,s) & \nabla_s\mathcal A_s(x,s)
    \end{bmatrix}.
\end{equation*}

Based on the above notation, we consider the following constraint
dissolving problem associated with \eqref{Eq:P}
\begin{equation}\tag{Q}\label{Eq:Q}
 \min_{x\in D,\ s\in C}
 \;h_\beta(x,s)
 :=f(\mathcal A_x(x,s))
   +\frac{\beta}{2}\|G(x)-s\|^2,
\end{equation}
where $\beta>0$ is a penalty parameter. 

In order to solve \eqref{Eq:P} via \eqref{Eq:Q}, it seems to be necessary to derive the equivalence of (M-)stationary points between \eqref{Eq:P} and \eqref{Eq:Q}. To this end, we first introduce the some notation for convenience and then present the notion of stationarity
for \eqref{Eq:Q}.

Let us define
\begin{equation}\label{eq:g-components}
 \begin{aligned}
 &g_x(x,s) :=\nabla_x\mathcal A_x(x,s)^\top
                                  \nabla f(\mathcal A_x(x,s)),\\
 &g_s(x,s) :=\nabla_s\mathcal A_x(x,s)^\top
                                  \nabla f(\mathcal A_x(x,s)),\\
 & g(x,s) :=\begin{pmatrix}g_x(x,s)\\g_s(x,s)\end{pmatrix}
       =\nabla\mathcal A(x,s)^\top
                     \begin{pmatrix}\nabla f(\mathcal A_x(x,s))\\0_{m\times 1}\end{pmatrix}.
 \end{aligned}
\end{equation}
The chain rule yields the two gradients
\begin{equation}\label{eq:penalty-gradients}
 \begin{aligned}
 \nabla_x h_\beta(x,s)&=g_x(x,s)+\beta \nabla G(x)^\top (G(x)-s),\\
 \nabla_s h_\beta(x,s)&=g_s(x,s)-\beta(G(x)-s).
 \end{aligned}
\end{equation}
We denote their column concatenation by $\nabla h_\beta(x,s)$.

\begin{definition}\label{def:penalty-stationarity}
A pair $(\bar x,\bar s)\in K$ is M-stationary for \eqref{Eq:Q} if
\begin{equation}\label{eq:m-penalty}
 \begin{aligned}
 0 \in\nabla_x h_\beta(\bar x,\bar s)+\mathcal N_D^{\lim}(\bar x)  \ {\rm and} \ 
 0\in\nabla_s h_\beta(\bar x,\bar s)+\mathcal N_C(\bar s).
 \end{aligned}
\end{equation}
\end{definition}

\begin{theorem}
\label{thm:penalty-m-equivalence}
Let Assumption~\ref{ass:P} hold and $\mathcal A$ be a
constraint dissolving mapping satisfying \Cref{def:am-regularity}. For any $\beta>0$ and
$\bar x\in\mathcal F$, one has
\begin{equation*}
    \bar x \text{ is M-stationary for \eqref{Eq:P}} \Longleftrightarrow \ (\bar x, G(\bar x))  \text{ is M-stationary for \eqref{Eq:Q}}.
\end{equation*}
\end{theorem}
\begin{proof}
Since $\bar x\in\mathcal F$, we have $(\bar x,G(\bar x))\in\mathcal F_s$, and
thus $\mathcal A_x(\bar x,G(\bar x))=\bar x$. Hence, by
\Cref{def:penalty-stationarity}, $(\bar x, G(\bar x))$ is an M-stationary
point of \eqref{Eq:Q} if and only if
\begin{equation*}
    0\in \nabla \mathcal A (\bar x, G(\bar x))^\top
    \begin{pmatrix}\nabla f(\bar x)\\ 0_m\end{pmatrix}
    + \mathcal N_K^{\lim}(\bar x, G(\bar x)).
\end{equation*}
By \Cref{thm:stationarity-equivalence}, this is equivalent to
\begin{equation}\label{eq:pf-range}
    0\in \begin{pmatrix}\nabla f(\bar x)\\ 0_m\end{pmatrix}
    + \operatorname{Range}\begin{pmatrix}\nabla G(\bar x)^\top\\ -I_{m}\end{pmatrix}
    + \begin{pmatrix} \mathcal N_D^{\text{lim}}(\bar x)\\ \mathcal N_C (G(\bar x))\end{pmatrix}.
\end{equation} 
On the other hand, $\bar x$ is an M-stationary point of \eqref{Eq:P} if and
only if there exists some $\lambda\in \mathbb R^m$ such that
\begin{equation}\label{eq:pf-lambda}
    0\in \begin{pmatrix}\nabla f(\bar x)\\ 0_m\end{pmatrix}
    +\begin{pmatrix}\nabla G(\bar x)^\top\\ -I_{m}\end{pmatrix}\lambda
    + \begin{pmatrix} \mathcal N_D^{\text{lim}}(\bar x)\\ \mathcal N_C (G(\bar x))\end{pmatrix},
\end{equation}
where the second block of \eqref{eq:pf-lambda} reads exactly
$\lambda\in\mathcal N_C(G(\bar x))$. Since
\begin{equation*}
    \operatorname{Range}\begin{pmatrix}\nabla G(\bar x)^\top\\ -I_m\end{pmatrix}
    =\left\{\begin{pmatrix}\nabla G(\bar x)^\top\\ -I_m\end{pmatrix}\xi:\
    \xi\in\mathbb R^m\right\},
\end{equation*}
\eqref{eq:pf-range} holds if and only if there exists
$\lambda\in\mathbb R^m$ such that \eqref{eq:pf-lambda} holds. Therefore,
$\bar x$ is an M-stationary point of \eqref{Eq:P} if and only if
$(\bar x,G(\bar x))$ is an M-stationary point of \eqref{Eq:Q}.
\end{proof}

\section{Algorithm}\label{Sec:alg}
Leveraging the equivalence of M-stationarity between \eqref{Eq:P} and \eqref{Eq:Ps}, this section develops a penalty-based constraint dissolving method for solving \eqref{Eq:Ps} and establishes its corresponding convergence analysis. The overall algorithm is presented in \Cref{Alg:CDP}.

\begin{algorithm}[H]\caption{A constraint dissolving method with adaptive penalty}
	\label{Alg:CDP}
	\begin{algorithmic}[1]
		\REQUIRE $\beta_0>0$, $\eta \in (0,1)$, $\tau>1$, $x^0\in D$, and $s^0=G(x^0)$.
            \STATE Set $k:=0$.
		\WHILE{$x^k$ does not satisfy a termination criterion}
        \STATE \label{Step:subpro}Compute an approximate M-stationary point $(x^{k+1}, s^{k+1})$ of the subproblem
        \begin{equation}\label{Eq:subph}
\begin{aligned}
  \min_{x,s} \ & h_{\beta_k}(x,s) \\ 
  \text{s.t.} \ & s\in C, x\in D,
\end{aligned}
\end{equation}
i.e., there exist
 $\varepsilon_x^{k+1}\in\mathbb R^n$ and
 $\varepsilon_s^{k+1}\in\mathbb R^m$ such that
\begin{equation}\label{eq:inner-stationarity}
 \begin{aligned}
 \varepsilon_x^{k+1}&\in\nabla_x h_{\beta_k}(x^{k+1},s^{k+1})
       +\mathcal N_D^{\lim}(x^{k+1}),\\
 \varepsilon_s^{k+1}&\in\nabla_s h_{\beta_k}(x^{k+1},s^{k+1})
       +\mathcal N_C(s^{k+1}).
 \end{aligned}
\end{equation}
        \IF{$k=0$ \OR $\|G(x^{k+1})-s^{k+1}\|\le\eta\|G(x^{k})-s^{k}\|$}
            \STATE $\beta_{k+1}:=\beta_k$
        \ELSE[\smash{$\left.\rule[-2.3\baselineskip]{0pt}{5\baselineskip}\right\}$}~Update $\beta_{k+1}$]
            \STATE $\beta_{k+1}:=\tau\beta_k$
        \ENDIF 
        \STATE Set $k:=k+1$.
        \ENDWHILE
		\RETURN $x^k$
	\end{algorithmic}
\end{algorithm}

\begin{remark}
(i) Step \ref{Step:subpro} requires the subproblem \eqref{Eq:subph} to be solved inexactly, for which we employ projected gradient descent method proposed in \cite[Algorithm~3.1]{jia2023augmented} whose either last iterate or any accumulation point is M-stationary.
\\ 
(ii) Regarding the penalty parameter $\beta_k$, it is updated in a nondecreasing manner according to the infeasibility measurement with respect to \eqref{Eq:Q}. This is in contrast to \cite{xiao2025exact}, where a fixed (large enough) penalty parameter is adopted, and \cite{jia2026constraint}, where the penalty parameter is allowed to vary dynamically but without any update rule.
\end{remark}

We first present the feasibility analysis of the accumulation points generated by \Cref{Alg:CDP}.

\begin{proposition}
\label{prop:bounded-values}
Let \Cref{ass:P} hold, the mapping $\mathcal A$ satisfy \Cref{def:cdm}, and $\{(x^k, s^k)\}$ be a sequence generated by \Cref{Alg:CDP}. Let $(\bar x, \bar s)$ be any accumulation point of $\{(x^k, s^k)\}$.
If one of the following conditions holds:
\begin{enumerate}[label=(\roman*)]
    \item \label{Item:bounded}$\{\beta_k\}$ is bounded, or
    \item \label{Item: B}there exists some $B\in\mathbb R$ such that $h_{\beta_k}(x^{k+1},s^{k+1})\leq B$ holds for all $k$,
\end{enumerate}
then $(\bar x, \bar s)$ is feasible for \eqref{Eq:Ps}  and $\bar x$ is feasible for \eqref{Eq:P}.
\end{proposition}
\begin{proof}
Let us choose an infinite index set $\mathcal K$ such that
$(x^{k},s^{k})\to_{\mathcal K}(\bar x, \bar s)$. By \eqref{Eq:subph} and the closedness of $C$ and $D$, we have $(\bar x, \bar s)\in D \times C$.

\ref{Item:bounded} If $\{\beta_k\}$ is bounded, by the update rule of $\{\beta_k\}$, we have $\|G(x^{k})-s^k\| \to 0$. The continuity of $G$ implies that $G(\bar x)=\bar s$. It means that $(\bar x, \bar s)$ is feasible for \eqref{Eq:Ps}. Furthermore, it implies that $\bar x\in D$ and $G(\bar x)=\bar s \in C$, which means that $\bar x$ is feasible for \eqref{Eq:P}. 

\ref{Item: B} 
By the definition of $h_{\beta_k}$ in \eqref{Eq:Q}, we have
\begin{equation}\label{Eq:h_beta_implies}
 \begin{aligned}
 0\leq\|G(x^k)-s^k\|^2
 &=\frac{2\bigl(h_{\beta_{k-1}}(x^{k},s^{k})
                     -f(\mathcal A_x(x^{k},s^{k}))\bigr)}{\beta_{k-1}}\\
 &\leq\frac{2\bigl(B-f(\mathcal A_x(x^{k},s^{k}))\bigr)}{\beta_{k-1}}, \ \forall \ k\geq 1.
 \end{aligned}
\end{equation}
Since $(x^k, s^k)\to_{\K} (\bar x, \bar s) \in D \times C$, the continuity of $f$ and $\mathcal A_x$ gives
\[
 f(\mathcal A_x(x^{k},s^{k}))\to_{\mathcal K}
 f(\mathcal A_x(\bar x,\bar s))=f(\bar x).
\]
Thus, taking the limit $k\to_{\K} \infty$ in \eqref{Eq:h_beta_implies} yields that $\|G(x^k)-s^k\| \to_{\K} 0$. The continuity of $G$ further implies that $G(\bar x)=\bar s$. Therefore, by the same analysis as the case that $\{\beta_k\}$ is bounded, we conclude that $(\bar x, \bar s)$ is feasible for \eqref{Eq:Ps} and $\bar x$ is feasible for \eqref{Eq:P}.
\end{proof}
 
\Cref{prop:bounded-values} provides the sufficient conditions for the feasibility of accumulation points. As illustrated in \cite{jia2023augmented,jia2026constraint}, \ref{Item: B} holds, for example, if $(x^{k+1}, s^{k+1})$ is an (approximate) global minimizer of the subproblem \eqref{Eq:subph}. This, however, requires more restrictive assumption on model, for instance, $h_{\beta}(x, s)$ is level bounded on $K$ \cite[Example~1.11]{rockafellar1998variational}. For the specific constraint dissolving mapping $\mathcal A(x,s)$ given by \eqref{eq:cdm-def}, every accumulation point is verified to be feasible under some mild conditions, see \Cref{Lem:feasibility_accumulation}, without requiring \ref{Item: B}. However, under a general constraint dissolving mapping $\mathcal A(x,s)$, we can only obtain that every accumulation point is an M-stationary point of the constraint violation problem.

\begin{proposition}
\label{prop:violation-stationarity}
Let \Cref{ass:P} hold, the mapping $\mathcal A$ satisfy \Cref{def:cdm}, and $\{(x^k, s^k)\}$ be a sequence generated by \Cref{Alg:CDP}. Let $(\bar x, \bar s)$ be any accumulation point of $\{(x^k, s^k)\}$. If  $\{(\varepsilon_x^k, \varepsilon_s^k)\}$ is bounded, then
$(\bar x, \bar s)$ is M-stationary of
\begin{equation}\label{eq:violation-problem}
 \min_{x\in D,\ s\in C}\ \frac12\|G(x)-s\|^2.
\end{equation}
\end{proposition}
\begin{proof}
As in the proof of \Cref{prop:bounded-values}, we take $\mathcal K$ is an
infinite index set such that $(x^k, s^k)\to_{\mathcal K}(\bar x, \bar s)$. 
Let us first consider the case $\{\beta_k\}$ is bounded, by \Cref{prop:bounded-values}, we have $G(\bar x)=\bar s$, which means that $(\bar x, \bar s)$ is a global minimizer of \eqref{eq:violation-problem} and hence its M-stationary point. It remains to consider the situation $\beta_k\to \infty$. 
By \eqref{eq:inner-stationarity}, one has for all $k\geq 1$, 
\[
 \begin{aligned}
 \frac{\varepsilon_x^{k}-g_x(x^{k},s^{k})}{\beta_{k-1}}
   &\in \nabla G(x^{k})^\top \left( G(x^k)-s^k \right)
                              +\mathcal N_D^{\lim}(x^{k}),\\
 \frac{\varepsilon_s^{k}-g_s(x^{k},s^{k})}{\beta_{k-1}}
   &\in-\left( G(x^k)-s^k \right)+\mathcal N_C(s^{k}). 
 \end{aligned} 
\]
By the boundedness of $\{(\varepsilon_x^k, \varepsilon_s^k)\}$, continuity of $\nabla G$, $G$, $\nabla \mathcal A$ and $\nabla f$, and the fact that $\mathcal N_D^{\lim}(\cdot)$ is robust, taking the limit $k\to_{\K} \infty$ implies that
\[ 
 0 \in \nabla G(\bar x)^\top(G(\bar x)-\bar s)
                                  +\mathcal N_D^{\lim}(\bar x), \ {\rm and} \
 0 \in-(G(\bar x)-\bar s)+\mathcal N_C(\bar s). 
\]
It means that $(\bar x, \bar s)$ is an M-stationary point of \eqref{eq:violation-problem}.
\end{proof}

\begin{theorem}\label{thm:main-convergence}
Let Assumption \ref{ass:P} hold, $\mathcal A$ satisfy \Cref{def:cdm}, and $\{(x^k, s^k)\}$ be a sequence generated by \Cref{Alg:CDP} with
$(\varepsilon_x^{k+1},\varepsilon_s^{k+1})\to(0,0)$. Let $(\bar x, \bar s)$ be any accumulation point of $\{(x^k, s^k)\}$ feasible for \eqref{Eq:Ps}. Then $\bar x$ is an M-stationary point for \eqref{Eq:P}.
\end{theorem}
\begin{proof}
    As in the proof of \Cref{prop:bounded-values}, we take $\mathcal K$ as an
infinite index set such that $(x^k, s^k)\to_{\mathcal K}(\bar x, \bar s)$. By the feasibility of $(\bar x, \bar s)$, we have $\bar s = G(\bar x)$.

By employing \eqref{eq:penalty-gradients}, \eqref{eq:inner-stationarity} reads
\[
 \begin{aligned}
 \varepsilon_x^{k}
 &\in g_x(x^{k},s^{k})
      +\beta_{k-1}\nabla G(x^{k})^\top H(x^{k},s^{k})
      +\mathcal N_D^{\lim}(x^{k}),\\
 \varepsilon_s^{k}
 &\in g_s(x^{k},s^{k})
      -\beta_{k-1} H(x^{k},s^{k})
      +\mathcal N_C(s^{k}), \  \forall \ k\geq 1.
 \end{aligned}
\]
It implies that
\[
 \begin{pmatrix}
  \varepsilon_x^{k}\\
  \varepsilon_s^{k}
 \end{pmatrix}
 -g(x^{k},s^{k})
 \in \operatorname{Range} \begin{pmatrix}
     \nabla G(x^k)^\top \\
     -I_m
 \end{pmatrix} + \mathcal N_K^{\text{lim}} (x^k, s^k),  \  \forall \ k\geq 1.\] 
In view of \eqref{eq:g-components}, \Cref{ass:P}, and \Cref{def:cdm}, \Cref{lem:multiplier-closedness} gives 
\begin{equation}\label{eq:limit-transformed-stationarity}
 0\in \nabla \mathcal A(\bar x, \bar s)^\top \begin{pmatrix}
     \nabla f (\bar x) \\
     0_{m}
 \end{pmatrix}
     +\operatorname{Range} \begin{pmatrix}
     \nabla G(\bar x)^\top \\
     -I_m
 \end{pmatrix}
     +\mathcal N_K^{\lim}(\bar x,\bar s).
\end{equation}
Let us recall that $\bar s = G(\bar x)$, it follows from \Cref{thm:penalty-m-equivalence} that $\bar x$ is an M-stationary point for \eqref{Eq:P}.
\end{proof}

\section{Implementation of constraint dissolving mappings for nonconvex sets}
\label{sec:realization}

This section focuses on the construction of constraint dissolving mappings $\mathcal A(x,s)$ satisfying \Cref{def:cdm}. To this end, we first present the motivation based on the idea from Fletcher \cite{fletcher1970class}.

Let us consider the following optimization problem
\begin{equation}\label{Eq:problemsz}
\begin{aligned}
 \min \ &f(w) \\
 \operatorname{s.t.}\ &c(w)=0, \ w\in K,
 \end{aligned}
\end{equation}
where $f$ and $c$ are continuously differentiable, and $K$ is nonempty, closed, and possibly nonconvex. For brevity, we do not specify the dimensions explicitly. All vectors, matrices and mappings are assumed to have compatible dimensions.

In order to solve \eqref{Eq:problemsz}, we alternatively solve
\begin{equation*}
\begin{aligned}
    \min \ & f(w) - c(w)^\top \lambda_w +\frac{\beta}{2}\|c(w)\|^2\\
    \operatorname{s.t.} \ & w\in K,
\end{aligned}
\end{equation*}
where $\lambda_w$ is a Lagrangian multiplier. The objective function takes the form of Fletcher's penalty function \cite{fletcher1970class}. 

Moreover, under some constraint qualification, the optimality condition of \eqref{Eq:problemsz} reads 
$0\in \nabla f(w)- \nabla c(w)^\top \lambda_w +\mathcal N_K^{\text{lim}}(w)$. 
 Then there exists some $\nu_w \in \mathcal N_K^{\text{lim}}(w)$ such that
 \begin{equation}\label{Eq:opt_z}
     0=\nabla f(w)- \nabla c(w)^\top \lambda_w + \nu_w.
 \end{equation}
 In fact, the optimal $\lambda_w$ depends on unknown normal
vector $\nu_w$ and hence cannot be evaluated. To remove the dependence on $\nu_w$, we assume that there exists a symmetric positive (semi)definite matrix $Q_w$ such that $\operatorname{Null}(Q_w)=\operatorname{Span}(\mathcal N_K^{\text{lim}}(w))$ (see \Cref{ass:projective}  \ref{item: Q} for more general setting). Then we have $Q_w \nu_w=0$. In this case, multiplying \eqref{Eq:opt_z} by $\nabla c(w)Q_w$ yields that
\begin{equation*}
    0= \nabla c(w)Q_w \nabla f(w)- \nabla c(w)Q_w \nabla c(w)^\top \lambda_w.
\end{equation*}
Under the condition that $\nabla c(w)Q_w \nabla c(w)^\top$ is invertible (a constraint qualification, like \Cref{Ass:A_CQ}, seems to be necessary to ensure this assumption, see \Cref{lem:invertibility}), one has
\begin{equation}\label{Eq:lambda_opt}
    \lambda_w= \left( \nabla c(w)Q_w\nabla c(w)^\top\right)^{-1} \nabla c(w)Q_w\nabla f(w).
\end{equation}

On the other hand, by employing the constraint dissolving method to solve \eqref{Eq:problemsz}, one has to solve
\begin{equation*}
\begin{aligned}
    \min \ &f(\mathcal A(w)) + \frac{\beta}{2}\|c(w)\|^2\\
    \operatorname{s.t.}\ & w\in K,
    \end{aligned}
\end{equation*}
where $\mathcal A$ is a constraint dissolving mapping. From the perspective that constraint dissolving method incorporates the multiplier correlation term from Fletcher's penalty function in \cite{fletcher1970class} implicitly through the constraint dissolving mapping $\mathcal A$, then
$f(\mathcal A(w)) \approx f(w) - c(w)^\top \lambda_w$. 
Since $f$ is continuously differentiable, one has 
$f(\mathcal A(w)) \approx f(w) + \left< \nabla f(w), \mathcal A(w)-w \right>$. 
By \eqref{Eq:lambda_opt}, it is reasonable to ask
\begin{equation*}
  \left< \nabla f(w), \mathcal A(w)-w \right> = -c(w)^\top \left( \nabla c(w)Q_w\nabla c(w)^\top\right)^{-1} \nabla c(w)Q_w\nabla f(w).  
\end{equation*}
Requiring this identity for every smooth objective $f$, i.e., for arbitrary
$\nabla f(w)$, leads to the candidate $\mathcal A(w)$ as follows
\begin{equation*}
    \mathcal A(w) := w - Q_w \nabla c(w)^ \top \left( \nabla c(w)Q_w\nabla c(w)^\top\right)^{-1} c(w).
\end{equation*}

Motivated by the above analysis, we impose the following assumption.

\begin{assumption}\label{ass:projective} 
\begin{enumerate}[label=(\alph*)]
 \item $G$ is twice continuously differentiable.
 \item \label{item: Q}There exist continuously differentiable mappings $\mathcal Q_D$ and $\mathcal Q_C$ defined on some open set containing $D$ and $C$, respectively, satisfying both $\mathcal Q_D(x)$ and $\mathcal Q_C(s)$ for all $x\in D$ and $s\in C$ are symmetric positive semidefinite such that 
       \begin{equation*}
       \begin{aligned}
        \operatorname{Null}\left(\mathcal Q_D(x)\right)&=\operatorname{Span}\left(\mathcal N_D^{\lim}(x)\right),&&x\in D,\\
        \operatorname{Null}\left(\mathcal Q_C(s)\right)&=\operatorname{Span}\left(\mathcal N_C(s)\right),&&s\in C.
       \end{aligned}
       \end{equation*}
\end{enumerate}
\end{assumption}

For several widely-used classes of $D$, including disjunctive, sparsity, low-rank-, low-rank positive semidefinite constraints, suitable choices of $\mathcal Q_D$ has been summarized in \cite[Table~1]{jia2026constraint}. Moreover, $C$ is usually given by a box-constraint, the corresponding choice of $\mathcal Q_C$ can be found in \cite[Table~2]{xiao2025exact}. We refer
interested readers to these references for more details.

For a fixed $\alpha > 0$, let us define
\begin{equation}\label{eq:block-matrices}
 \begin{aligned}
 &\mathcal Q(x,s) :=\begin{pmatrix}\mathcal Q_D(x)&0\\0&\mathcal Q_C(s)\end{pmatrix},\\
 &M_0(x,s) :=\nabla H(x,s)\mathcal Q(x,s)\nabla H(x,s)^\top =\nabla G(x)\mathcal Q_D(x)\nabla G(x)^\top+\mathcal Q_C(s),\\
 &M(x,s) :=M_0(x,s)+\alpha\|H(x,s)\|^2I_m.
 \end{aligned}
\end{equation}
Based on the above notation, for $(x,s)\in K$, we define
\begin{equation}\label{eq:cdm-def}
 \mathcal A(x,s):=(x,s)-\mathcal Q(x,s)\nabla H(x,s)^\top M(x,s)^{-1}H(x,s),
\end{equation}
Note that $M_0(x,s)$ is positive definite on $K$ under some constraint qualification (see \Cref{lem:invertibility}); however, it may be nearly singular.
The regularization term $\alpha\|H(x,s)\|^2I_m$ alleviates this issue. Indeed, for every $(x,s)\in K\setminus\mathcal F_s$,
\[
  \bigl\|M(x,s)^{-1}H(x,s)\bigr\|
  \le\frac{\|H(x,s)\|}{\lambda_{\min}(M_0(x,s))+\alpha\|H(x,s)\|^2}
  \le\frac{1}{\alpha\|H(x,s)\|},
\]
independent of $\lambda_{\min}(M_0(x,s))$. However, this term vanishes
on $\mathcal F_s$ regardless of the value of $\alpha$, so that $M(x,s)=M_0(x,s)$
there.

We next present a counterexample showing that the AM-regularity, introduced in \Cref{def:am-regularity} and used in the preceding sections, is not sufficient to guarantee the positive definiteness of $M_0(x,s)$ and $M(x,s)$ on $K$.

\begin{example}\label{ex:AM-not-enough}
Let $n=1$, $m=2$, $f\equiv0$, $D=\mathbb R$, $C=\mathbb R^2_-$, and $G(x)=(x,-x)^\top$. Then
$\mathcal F=\{0\}$ and $\mathcal F_s=\{(0,(0,0))\}$.  
The point $\bar x=0$ is AM-regular. Indeed, since
$\mathcal N_D^{\text{lim}}(x)\equiv\{0\}$, we have 
$\mathcal M(x,z)=\{\lambda_1-\lambda_2:\ \lambda \in \mathcal N_C(G(x)-z)\}$. 
By $\mathcal N_C((0,0))=\mathbb R_{+}^2$, we have $\mathcal M(0,0)=\mathbb R$, then
$\limsup_{x\to0,\,z\to0}\mathcal M(x,z)\subseteq\mathcal M(0,0)$ holds
trivially. 
On the other hand, to ensure \Cref{ass:projective} \ref{item: Q}, we set $\mathcal Q_C(0)=0$ and choose $\mathcal Q_D(0)=1$. By $\nabla G(x)\equiv(1,-1)^\top$ and
$H(0,(0,0))=0$, for every $\alpha\geq0$, we have
\[
 M(0,(0,0))=M_0(0,(0,0))=\nabla G(0)\nabla G(0)^\top
 =\begin{pmatrix}1&-1\\-1&1\end{pmatrix},
\]
which is singular.
Therefore, AM-regularity is not sufficient to ensure that $M(x,s)$ is positive
definite on $\mathcal F_s$.
\end{example}

Therefore, a stronger constraint qualification is needed to ensure $M_0(x,s)\succ 0$ and thus $M(x,s) \succ 0$ for every $(x,s)\in K$. To this end, we impose the following condition, which generalizes the Robinson nondegeneracy condition \cite[Definition~2.1, (2.3)]{robinson2003constraint}. In particular, the two conditions coincide when $C$ and $D$ are two polyhedral convex sets.

\begin{assumption}\label{Ass:A_CQ} 
 At any $(x,s)\in K$,  $d=0$ is the only vector
       $d\in\mathbb R^m$ satisfying
       \begin{equation*} 
        \nabla G(x)^\top d \in\operatorname{Span}\bigl(\mathcal N_D^{\lim}(x)\bigr) \ {\rm and} \
        d \in\operatorname{Span}\bigl(\mathcal N_C(s)\bigr).  
    \end{equation*} 
\end{assumption} 

\begin{lemma}\label{lem:invertibility}
Suppose that \Cref{ass:projective} and \Cref{Ass:A_CQ} hold, $M_0(x,s)\succ0$ and thus $M(x,s)\succ0$ at every $(x,s)\in K$.
Moreover, every $\bar x\in\mathcal F$ is AM-regular.
\end{lemma}
\begin{proof}
For any $(x,s)\in K$, as well as $d\in\mathbb R^m$,
\[
 d^\top M_0(x,s)d
  =(\nabla G(x)^\top d)^\top \mathcal Q_D(x)(\nabla G(x)^\top d)
       +d^\top \mathcal Q_C(s)d\geq0.
\]
Hence $M_0(x,s) \succeq 0$ for all $(x, s)\in K$.

Now suppose $d^\top M_0(x,s)d =0$. Then the positive semidefiniteness of $\mathcal Q_D(x)$ and $\mathcal Q_C(s)$ implies
\[
 \nabla G(x)^\top d\in\operatorname{Null}(\mathcal Q_D(x)),
 \  d\in\operatorname{Null}(\mathcal Q_C(s)).
\]
By \Cref{ass:projective} \ref{item: Q} and \Cref{Ass:A_CQ}, one has
$d=0$. Thus $M_0(x,s)\succ0$ for all $(x, s)\in K$. Since $\alpha\|H(x,s)\|^2I_m \succeq 0$ by $\alpha > 0$, $M(x,s) \succ 0$ for all $(x,s)\in K$.

For any $\bar x\in \mathcal F$, fix some $\xi \in\operatorname{limsup}_{x\to \bar x, z\to 0} \mathcal M(x, z)$. Then we find sequences $\{x^k\}\subset \mathbb R^n$, $\{\xi^k\} \subset \mathbb R^n$, and $\{z^k\}\subset \mathbb R^m$ such that $x^k\to_D \bar x$, $\xi^k\to \xi$, $z^k \to 0$, and $\xi^k\in \mathcal M(x^k, z^k)$. In particular, there exist sequences $\{\lambda^k\}\subset \mathbb R^m$ and $\{u^k\}\subset \mathbb R^n$ satisfying $\lambda^k\in\mathcal N_C(G(x^k)-z^k)$ and
 $u^k\in\mathcal N_D^{\lim}(x^k)$ such that
\begin{equation}\label{Eq:xi_AMreg}
\xi^k=\nabla G(x^k)^\top\lambda^k+u^k, \  \forall \ k\in \mathbb N.
\end{equation} 
If $\{(\lambda^k, u^k)\}$ is unbounded, without loss of generality, we assume
$\|(\lambda^k, u^k)\|\to\infty$.  Then it is reasonable to set
$(\lambda^k, u^k)/\|(\lambda^k, u^k)\|\to(\lambda, u) \neq (0,0)$.
Dividing \eqref{Eq:xi_AMreg} by $\|(\lambda^k, u^k)\|$, by the continuity of $G$ and the outer semicontinuity of the limiting normal cone, we have 
$\lambda\in\mathcal N_C(G(\bar x))$ and $-\nabla G(\bar x)^\top\lambda= u\in\mathcal N_D^{\lim}(\bar x)$. 
Then \Cref{Ass:A_CQ} implies $\lambda =0$, and thus $u= -\nabla G(\bar x)^\top\lambda=0$, which contradicts $(\lambda, u) \neq (0,0)$. Thus, the sequence $\{(\lambda^k, u^k)\}$ is bounded, and therefore, possesses a convergent subsequence with limit $(\hat \lambda, \hat u)$. By the continuity of $\nabla G$, taking the limit in \eqref{Eq:xi_AMreg} along this subsequence yields that
\begin{equation}\label{Eq:xi_AM_reg_right}
    \xi = \nabla G(\bar x)^ \top \hat \lambda + \hat u.
\end{equation}
Moreover, by the continuity of $G$ and $\nabla G$, $z^k\to 0$, $x^k\to_D \bar x$, and the outer semicontinuity of the limiting normal cone, we have $\hat \lambda \in \mathcal N_C (G(\bar x))$ and $\hat u \in \mathcal N_D^{\text{lim}}(\bar x)$. Thus, \eqref{Eq:xi_AM_reg_right} implies that 
$\xi \in \nabla G(\bar x)^ \top \mathcal N_C (G(\bar x)) + \mathcal N_D^{\text{lim}}(x) =\mathcal M(\bar x, 0)$. 
This shows that $\bar x$ is AM-regular.
\end{proof}

Note that \Cref{lem:invertibility} ensures the well-definedness of $\mathcal A$
in \eqref{eq:cdm-def}. We next verify that it satisfies \Cref{def:cdm}.
For brevity, we denote
\begin{equation}\label{eq:correction-matrix}
 J(x,s):=\mathcal Q(x,s)\nabla H(x,s)^\top M(x,s)^{-1}.
\end{equation}

\begin{proposition}\label{prop:block-cdm}
Suppose that \Cref{ass:P}, \Cref{ass:projective}, and \Cref{Ass:A_CQ} hold, the mapping
\eqref{eq:cdm-def} is a constraint dissolving mapping, i.e., it satisfies \Cref{def:cdm}.
\end{proposition}
\begin{proof}
    Since the proof is similar to that in \cite[Proposition~4.3]{jia2026constraint}, we present it in Appendix \ref{app:proof-block-cdm}.
\end{proof}
We next show that, with $\mathcal A$ specifically given by \eqref{eq:cdm-def}, every
accumulation point of the sequence generated by \Cref{Alg:CDP} is feasible for \eqref{Eq:Ps},
provided that the inexactness sequence $\{(\varepsilon_x^k,\varepsilon_s^k)\}$
is bounded.

\begin{lemma}\label{Lem:feasibility_accumulation}
    Suppose that \Cref{ass:P}, \Cref{ass:projective}, and \Cref{Ass:A_CQ} hold, the constraint dissolving mapping $\mathcal A(x,s)$ is defined by \eqref{eq:cdm-def}. Let $\{(x^k, s^k)\}$ be generated by \Cref{Alg:CDP} with $\{(\varepsilon_x^{k}, \varepsilon_s^k)\}$ being bounded, and $(\bar x, \bar s)$ be any accumulation point of $\{(x^k, s^k)\}$. Then $(\bar x, \bar s)$ is feasible for \eqref{Eq:Ps} and thus $\bar x$ is feasible for \eqref{Eq:P}.
\end{lemma}
\begin{proof}
By \Cref{prop:bounded-values}~\ref{Item:bounded}, it suffices to consider the case where $\{\beta_k\}$ is unbounded. 
Let us also take $\mathcal K$ as an infinite index set such that $(x^k, s^k) \to_{\mathcal K}(\bar x, \bar s)$, and denote $ w:=(x, s)$ and $\bar w:=(\bar x, \bar s)$ for convenience. By the closedness of $C$ and $D$, we have $\bar w\in K$. Then by \Cref{lem:invertibility}, we have $M_0(\bar w)\succ0$. Let us set
$\mu:=\lambda_{\min}(M_0(\bar w))>0$. In view of \Cref{ass:projective}, $M_0(\cdot)$, and thus $M(\cdot)$ defined in \eqref{eq:block-matrices} are continuous, then there exists $\rho_1>0$ such that
$
  M(w)\succeq M_0(w) \succeq \frac{\mu}{2}I_m$ on $B_1:=B_{\rho_1}(\bar w)$.
By \Cref{ass:P} and \Cref{ass:projective}, $J$ (defined in \eqref{eq:correction-matrix}), $DJ$, $\mathcal A$, $H$, and $\nabla f$ are continuous. Then, the compactness of $B_1$ implies that the constants
\begin{equation*}
    \begin{aligned}
         L_J&:=\max_{w\in B_1}\|J(w)\|, \quad
  L_{DJ}:=\max_{w\in B_1}\|DJ(w)\|, \\
  L_H&:= \max_{w\in B_1} \|H(w)\|, \quad
   L_f:=\max_{w\in B_1}\|\nabla f(\mathcal A_x(w))\|
    \end{aligned}
\end{equation*}
are finite.  
Recall that $w^k:=(x^k, s^k)\to_{\K} (\bar x, \bar s):=\bar w$. For all large enough $k\in \K$, one has $w^k\in B_{\rho}(\bar w)$. By \Cref{ass:projective} and \eqref{eq:correction-matrix}, one has
\begin{equation}\label{Eq:Jlimcone=Null}
    J(w^k)^\top\nu^k=M(w^k)^{-1}\nabla H(w^k)\mathcal Q(w^k)\nu^k=0, \ \forall \ \nu^k\in \mathcal N_K^{\text{lim}}(w^k), 
\end{equation}
and thus
\begin{equation}\label{eq:pf-normal}
  \big\langle\nu^k,\,J(w^k)H(w^k)\big\rangle
  =\big\langle J(w^k)^\top\nu^k,\,H(w^k)\big\rangle=0 .
\end{equation} 
Moreover, by \eqref{eq:block-matrices} and \eqref{eq:correction-matrix}, we have
\begin{equation}\label{eq:pf-identity}
\begin{aligned}
  \nabla H(w^k)J(w^k)
  &=\nabla H(w^k)\mathcal Q(w^k)\nabla H(w^k)^\top M(w^k)^{-1}
  =M_0(w^k)M(w^k)^{-1}\\
 & =I_m-\alpha\|H(w^k)\|^2M(w^k)^{-1}.
  \end{aligned}
\end{equation} 
By \eqref{eq:inner-stationarity}, there exists some $v^k\in \mathcal N_K^{\text{lim}}(w^k)$ such that $\varepsilon^k=\nabla h_{\beta_{k-1}}(w^k)+v^k$ with $\varepsilon^k:=(\varepsilon_x^k, \varepsilon_s^k)$ for each $k\geq 1$. Consequently, we have
\begin{equation*}
   0= \left< \nabla h_{\beta_{k-1}}(w^k)+v^k-\varepsilon^k, J(w^k)H(w^k) \right>, \  \forall \ k\geq 1.
\end{equation*}
Let us write $\xi^k:=J(w^k)H(w^k)$ for brevity, then \eqref{Eq:Jlimcone=Null} implies
\begin{equation}\label{Eq:innerprod_xi=Null}
0=\left< \nabla \mathcal A_x(w^k)^\top \nabla f(\mathcal A_x(w^k)), \xi^k \right> + {\beta_{k-1}} \left<\nabla H(w^k)^\top H(w^k), \xi^k \right>-\left<\varepsilon^k, \xi^k \right>
\end{equation}
for all large enough $k \in \K$. Moreover, we have $\|\xi^k\| \leq L_J \|H(w^k)\|$ for all large enough $k \in \K$. 
By \eqref{eq:pf-identity}, we have
\begin{equation*}
\begin{aligned}
    J(w^k) \nabla H(w^k) \xi^k &= J(w^k) \nabla H(w^k) J(w^k)H(w^k) \\
    &= J(w^k)H(w^k)-\alpha \|H(w^k)\|^2 J(w^k)M(w^k)^{-1} H(w^k)
    \end{aligned}
\end{equation*}
for all $k \geq 1$. Since $\mathcal A$ is defined by \eqref{eq:cdm-def}, then we have
\begin{equation*}
\begin{aligned}
    \|\nabla\mathcal A(w^k)\xi^k\|&= \left\|\xi^k -\big(\mathcal D J(w^k)[\xi^k]\big)H(w^k)-J(w^k) \nabla H(w^k) \xi^k\right\|\\
    & = \left\|\big(\mathcal D J(w^k)[\xi^k]\big)H(w^k)-\alpha \|H(w^k)\|^2 J(w^k)M(w^k)^{-1} H(w^k)\right\|\\
 & \le L_{DJ}L_J\|H(w^k)\|^2+ \frac{2\alpha L_J}{\mu}\|H(w^k)\|^3\\
&  \le \left( L_{DJ}L_J + \frac{2\alpha L_JL_H}{\mu} \right)\|H(w^k)\|^2,  \ \forall \ \text{large enough}\  k\in \K.  
  \end{aligned}
\end{equation*}
Then 
\begin{equation}\label{Eq:Uppbound_fAinnerp}
\left< \nabla \mathcal A_x(w^k)^\top \nabla f(\mathcal A_x(w^k)), \xi^k \right> \leq L_f \left( L_{DJ}L_J + \frac{2\alpha L_JL_H}{\mu} \right)\|H(w^k)\|^2
\end{equation}
for all large enough $k\in \K$.  
On the other hand, by the notation of $M_0$ and $M$ in \eqref{eq:block-matrices}, one has $$\lambda_{\operatorname{min}}\left(M_0(w^k)M(w^k)^{-1}\right) \geq \frac{\mu}{\mu+2\alpha L_H^2}.$$ 
Then, in view of \Cref{lem:invertibility} and \eqref{eq:pf-identity}, we have
\begin{equation}\label{Eq:Lowbound_H}
\begin{aligned}
    \left<\nabla H(w^k)^\top H(w^k), \xi^k \right>&=\left<H(w^k), \nabla H(w^k) J(w^k)H(w^k)\right>\\
    &= H(w^k)^\top M_0(w^k)M(w^k)^{-1} H(w^k)\\
    & \geq \frac{\mu}{\mu+2\alpha L_H^2} \|H(w^k)\|^2
    \end{aligned}
\end{equation}
for all large enough $k\in \K$. 
Since $\{\varepsilon^k\}$ is bounded, we let $E>0$ be the constant satisfying $\|\varepsilon^k\| \leq E$ for all sufficiently large $k\in \K$. Then by \eqref{Eq:innerprod_xi=Null}, \eqref{Eq:Uppbound_fAinnerp}, and \eqref{Eq:Lowbound_H}, one has
\begin{equation*}
    \beta_{k-1} \frac{\mu}{\mu+2\alpha L_H^2} \|H(w^k)\|^2 \leq L_f\left(L_{DJ}L_J + \frac{2\alpha L_JL_H}{\mu} \right)\|H(w^k)\|^2 + E L_J L_H,
\end{equation*}
and consequently,
\begin{equation*}
    \left(  \beta_{k-1} \frac{\mu}{\mu+2\alpha L_H^2} - L_fL_{DJ}L_J -L_f\frac{2\alpha L_JL_H}{\mu} \right) \|H(w^k)\|^2 \leq E L_J L_H
\end{equation*}
for all large enough $k\in \K$. Since $\beta_{k-1}=\beta_k/\tau \to \infty$, then $H(w^k)\to H(\bar w)=0$. This means that $(\bar x, \bar s)$ is feasible for \eqref{Eq:Ps}, and thus $\bar x$ is feasible for \eqref{Eq:P}.
\end{proof}

By \Cref{Lem:feasibility_accumulation} and \Cref{thm:main-convergence}, we directly have the following result.

\begin{theorem}
Let \Cref{ass:P}, \Cref{ass:projective}, and \Cref{Ass:A_CQ} hold, $\mathcal A$ be defined in \eqref{eq:cdm-def}, and $\{(x^k, s^k)\}$ be the sequence generated by \Cref{Alg:CDP} with
$(\varepsilon_x^{k+1},\varepsilon_s^{k+1})\to(0,0)$. Let $(\bar x, \bar s)$ be any accumulation point of $\{(x^k, s^k)\}$. Then $\bar x$ is an M-stationary point for \eqref{Eq:P}.
\end{theorem}

\section{Numerical Experiments}\label{Sec:Num}
In this section, we report numerical results for \Cref{Alg:CDP}, referred to as CD for short, in comparison with the safeguarded augmented Lagrangian method (ALM) \cite{jia2023augmented} and the inexact penalty decomposition (PD) method  \cite{kanzow2023inexact}. All experiments were carried out in MATLAB (R2025b) on a MacBook Air with an Apple M4 chip and 16 GB of memory.

CD and ALM solve their subproblems by the projected gradient method
\cite[Algorithm~3.1]{jia2023augmented}, where the parameters are the same as those in \cite[Section~6]{jia2023augmented}.
For CD, in turn, the parameters $\eta$ and $\tau$ in \Cref{Alg:CDP},
the parameter $\alpha$ of the constraint dissolving mapping, and the
mappings $\mathcal Q_D$ and $\mathcal Q_C$ depend on the problem and
will be specified below. PD replaces the constraint $x\in D$ by $x-y=0$ and $y\in D$ and
penalizes the constraints $x-y=0$ and $G(x)\in C$, so that its
subproblems consist of minimizing $q_\beta(x,y):=f(x)+\frac\beta2\left(\|x-y\|_2^2
+\operatorname{dist}_C(G(x))^2\right)$
with respect to $x\in\mathbb R^n$ and $y\in D$.
Following \cite[Algorithm~2]{kanzow2023inexact}, we solve these
subproblems by alternating a gradient step in $x$, whose step size is
determined by an Armijo line search,  
and the update $y=P_D(x)$.

In outer iteration $k\geq 0$, the inner stopping criterion of CD and
ALM is
\cite[(3.4)]{jia2023augmented} with tolerance
$\varepsilon_{\mathrm{inner},k}=10^{-4}/\sqrt{k+1}$, while the
alternating scheme of PD is terminated when
$q_{\beta_k}(x^k,y^k)-q_{\beta_k}(x^{k+1}, y^{k+1})\leq 10^{-5}$.
The outer stopping criteria are as follows: CD terminates when
$\|G(x^{k+1})-s^{k+1}\|_\infty\leq\varepsilon_{\mathrm{out}}$,
ALM when \cite[(4.5)]{jia2023augmented} is satisfied, and PD when $\|x^{k+1}-y^{k+1}\|_\infty
+\|G(x^{k+1})-P_C(G(x^{k+1}))\|_\infty\leq\varepsilon_{\mathrm{out}}$
\cite[Section~6]{kanzow2023inexact}. The same tolerance
$\varepsilon_{\mathrm{out}}$ is used for all three methods, with $\varepsilon_{\mathrm{out}}=10^{-4}$ for the problems in
\Cref{subsec:num-rosenbrock,subsubsec:num-maxcut} and $\varepsilon_{\mathrm{out}}=10^{-5}$
for the real-world portfolio problems in \Cref{sec:portfolio-models}.
In addition, we impose the upper bound $10^8$ on the
penalty parameter and a budget of $10^7$ evaluations of the subproblem
objective functions per run. 
Given a point $\widetilde x^0$, all three methods are initialized
from the same projected point $x^0\in P_D(\widetilde x^0)$. CD
additionally sets $s^0=P_C(G(x^0))$, and PD sets $y^0=x^0$.
Following \cite[Section~6]{jia2023augmented}, ALM and PD use the same initial penalty
parameter as therein,
while the initial penalty parameter of CD is specified individually in the following subsections.

\subsection{A complementarity-constrained Rosenbrock problem}
\label{subsec:num-rosenbrock}

We first consider the following Rosenbrock problem with interval and complementarity constraints:
\begin{equation}
	\label{eq:num-rosenbrock}
	\begin{aligned}
		&\min_{x\in\mathbb R^2}\quad
		   f(x):=(1-x_1)^2+100(x_2-x_1^2)^2,\\
		&\text{s.t.} \quad
		 \frac12\leq x_1+x_2\leq\frac32, \ x_1\geq0,\ x_2\geq0,\ x_1x_2=0,
	\end{aligned}
\end{equation}
which corresponds to $G(x)=x_1+x_2$, $C=\left[\frac12,\frac32\right]$,
and $D=\{x\in\mathbb R_+^2:x_1x_2=0\}$. Its global
minimizer is $x^\star=(1/2,0)^\top$ with $f(x^\star)=6.5$ and
local minimizer is $\widehat x=(0,1/2)^\top$ with $f(\widehat x)=26$.
For the constraint dissolving mapping, we choose
$\mathcal Q_D(x)=c\,\diag(x_1^2,x_2^2)$ and
$\mathcal Q_C(s)=\left(s-\frac12\right)\left(\frac32-s\right)$, and we test the values
$c\in\{0.01,0.1,1,10,100\}$ in order to illustrate the corresponding numerical influence.

All three methods are run from the same $1000$
starting points $\widetilde x^{\,i}=(a_i,a_i^2+0.1\xi_i)^\top$, where
 $a_i$ is equally spaced in $[-1,1]$ and $\xi_i$ are drawn
independently from the uniform distribution on $[0,1]$. In this problem, we set $\tau=10$, $\eta=0.5$, and $\alpha=10$ for CD, and we use the same initial penalty parameter $\beta_0$ for the three methods.

The results are summarized in \Cref{tab:num-rosenbrock}, which
reports the average numbers of outer iterations ($Outer$) and of
cumulative inner iterations ($Inner$).
$N_{\mathrm{global}}$ and $N_{\mathrm{local}}$ are the numbers of runs whose final
iterate has a distance of at most $10^{-3}$ to $x^\star$ and
$\widehat x$, respectively, and the average running time per run in
milliseconds ($Time$). 
For each runs, the final iterate is close to either $x^\star$ or
$\widehat x$ in the above sense.
However, CD finds the global minimizer more often: for every value
of $c$, it terminates near $x^\star$ in $813$ to $968$ of the $1000$
runs, whereas ALM and PD do so in $668$ and $800$ runs, respectively.
Moreover, CD requires fewer iterations. Its average number of outer iterations decreases from $4.050$ to $1.298$ as $c$ increases from $0.01$ to $100$, compared with $15.488$ for ALM and $104.795$ for PD, and its average number of inner iterations lies between $16.170$ and $33.002$, compared with $56.778$ for ALM and $3075.530$ for PD.
The number of inner iterations and the running time of CD are
smallest for the intermediate values $c=1$ and $c=10$, for which CD
requires $0.40$ and $0.39$~ms on average, respectively, compared with
$0.52$~ms for ALM and $37.86$~ms for PD. 

\begin{table}[t]
	\centering
	\caption{Averaged results of CD, ALM, and PD for solving \eqref{eq:num-rosenbrock}.}\label{tab:num-rosenbrock}
	\begin{tabular}{ccccccc}
        \toprule
		Alg.
		& $c$
		& $Outer$
		& $Inner$
		& $N_{\mathrm{global}}$
		& $N_{\mathrm{local}}$
		& $Time$\\
		\midrule
		CD & $0.01$ & $4.050$ & $32.154$ & $968$ & $32$ & $0.58$ \\
		CD & $0.1$ & $2.963$ & $22.569$ & $934$ & $66$ & $0.51$ \\
		CD & $1$ & $2.035$ & $16.170$ & $906$ & $94$ & $0.40$ \\
		CD & $10$ & $1.458$ & $16.681$ & $831$ & $169$ & $0.39$ \\
		CD & $100$ & $1.298$ & $33.002$ & $813$ & $187$ & $1.23$ \\
		\noalign{\vskip\aboverulesep}\hdashline\noalign{\vskip\belowrulesep}
		ALM & -- & $15.488$ & $56.778$ & $668$ & $332$ & $0.52$ \\
		PD & -- & $104.795$ & $3075.530$ & $800$ & $200$ & $37.86$ \\
\bottomrule
	\end{tabular}
\end{table}

\subsection{Sparse portfolio optimization}
\label{sec:portfolio-models}
This subsection considers two types of portfolio problems with cardinality constraints, i.e., the interval-budget model in \Cref{subsubsec:portfolio-interval} and the mean-variance model in \Cref{subsec:num-mv}, as the variants of the classical portfolio problem \cite{bertsimas2022scalable}. We denote
$\Sigma$ and $\mu$ as the covariance matrix and the vector of
expected returns of $n$ assets and $e$ as the all-one vector.
In both problems, we use $D=\{x\in\mathbb R_+^n:\|x\|_0\leq\kappa\}$
and $\mathcal Q_D(x)=c\,\diag(x_1^2,\ldots,x_n^2)$.
For all three methods, we choose the starting point as $x^0=P_D(\widetilde x^0)$, where $\widetilde x^0$ is a
solution of the respective problem without the cardinality constraint
$\|x\|_0\leq\kappa$, which is a convex quadratic program.
For both problems, PD fails to satisfy the outer stopping criterion, and its results are therefore not reported.

\subsubsection{Interval-budget model on benchmark instances}
\label{subsubsec:portfolio-interval}
We first consider the sparse portfolio problem with an interval budget
constraint
\begin{equation}\label{eq:portfolio-interval}
	\begin{aligned}
		\min_{x\in\mathbb R^n}\quad & f(x):=\frac12x^\top\Sigma x-\mu^\top x \\
		\text{s.t.}\quad &0.8\leq e^\top x\leq1, \  x\geq0,
		\  \|x\|_0\leq\kappa,
	\end{aligned}
\end{equation}
in which at least $80\%$ of the capital has to be invested and the
remainder may be held in cash.
In contrast to the usual budget constraint $e^\top x=1$, under which
the whole capital is invested, this yields a test problem in which $C$
is not a singleton.
More precisely, we set $G(x)=e^\top x-1$, $C=[-0.2,0]$, and
$\mathcal Q_C(s)=(s+0.2)(-s)$.
The data $\Sigma$ and $\mu$ are taken from the test problem collection
\cite{frangioni2007sdp}, from which we use the orl and pard instances
of dimension $n\in\{200,300,400\}$, $30$ ($20$ orl
and $10$ pard) instances for each dimension, together with $\kappa\in\{3,5,7\}$, resulting
in $270$ test problems.
For CD, we set $\tau=2$, $\eta=0.5$, $c=0.01$, $\alpha=0.001$, and
$\beta_0=300$.

\begin{figure}[tbp]
	\centering
	\includegraphics[width=\textwidth]{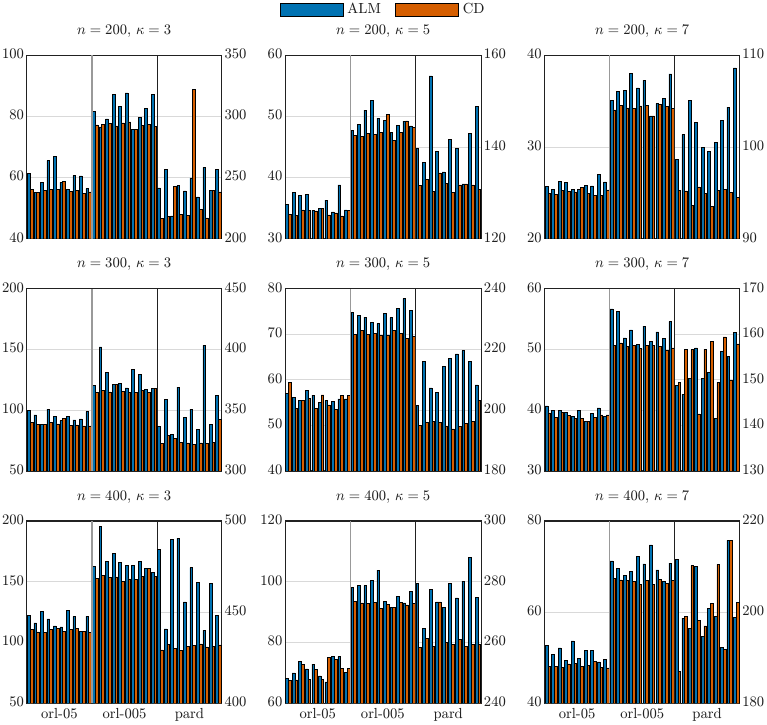}
	\caption{Final objective values of ALM and CD for solving
		\eqref{eq:portfolio-interval}. Each panel shows, for one pair
		$(n,\kappa)$, the $20$ orl instances (left part, left axis) and
		the $10$ pard instances (right part, right axis) of
		\cite{frangioni2007sdp}. Note that the gray vertical line separates the two
		families of orl instances.}
	\label{fig:num-portfolio}
\end{figure}

The results are summarized in \Cref{tab:num-portfolio}, which
reports, for each pair $(n,\kappa)$, the average numbers of outer
iterations ($Outer$) and of cumulative inner iterations ($Inner$), the
average running time in milliseconds ($Time$, minimum over three
repetitions), the average final objective function value ($\bar f$),
and the numbers of instances on which CD attains the lower objective
function value ($N_{\mathrm{CD}}$), ALM attains the lower one
($N_{\mathrm{ALM}}$), and the two methods tie ($N_{\mathrm{tie}}$).
In addition, Figure~\ref{fig:num-portfolio} shows the final objective
function values for every instance. 
Regarding the objective function values, CD attains the lower value
on $226$ of the $270$ instances and ALM on $30$, while the two methods
tie on the remaining $14$.
In each of the nine settings, CD attains the lower value on at least
$21$ of the $30$ instances and its average objective function value is
lower than that of ALM.
Regarding the numerical effort, CD requires fewer
outer iterations, fewer inner iterations, and less running time
than ALM in all nine settings. More precisely, its average number of inner iterations
is between $33\%$ and $60\%$ of that of ALM.
Over all settings, the average running times range from $3.3$ to
$6.6$~ms for CD and from $4.5$ to $10.7$~ms for ALM. 

\begin{table}[t]
	\centering
	\caption{Averaged results of CD and ALM for solving \eqref{eq:portfolio-interval}.
		}
	\label{tab:num-portfolio}
	\footnotesize\setlength{\tabcolsep}{3.5pt}
	\begin{tabular}{cc cccc cccc ccc}
		\toprule
		& & \multicolumn{4}{c}{CD} & \multicolumn{4}{c}{ALM} & & & \\
		\cmidrule(lr){3-6}\cmidrule(lr){7-10}
		$n$ & $\kappa$
		& $Outer$ & $Inner$ & $Time$ & $\bar f$
		& $Outer$ & $Inner$ & $Time$ & $\bar f$
		& $N_{\mathrm{CD}}$ & $N_{\mathrm{ALM}}$ & $N_{\mathrm{tie}}$\\
		\midrule
		$200$ & $3$ & $6.07$ & $138.8$ & $3.7$ & $122.943$ & $27.30$ & $335.0$ & $6.1$ & $128.637$ & $23$ & $5$ & $2$ \\
		$200$ & $5$ & $7.53$ & $124.3$ & $3.3$ & $71.241$ & $19.27$ & $259.5$ & $4.5$ & $75.488$ & $26$ & $1$ & $3$ \\
		$200$ & $7$ & $7.67$ & $125.4$ & $3.7$ & $51.399$ & $26.63$ & $376.0$ & $6.8$ & $54.770$ & $28$ & $2$ & $0$ \\
		$300$ & $3$ & $5.93$ & $129.0$ & $3.6$ & $176.978$ & $32.33$ & $395.9$ & $8.9$ & $191.376$ & $26$ & $1$ & $3$ \\
		$300$ & $5$ & $7.27$ & $151.1$ & $4.6$ & $107.272$ & $25.57$ & $354.2$ & $7.6$ & $114.332$ & $25$ & $4$ & $1$ \\
		$300$ & $7$ & $8.17$ & $167.7$ & $5.5$ & $81.081$ & $20.17$ & $278.8$ & $6.2$ & $81.603$ & $21$ & $8$ & $1$ \\
		$400$ & $3$ & $6.80$ & $154.5$ & $5.3$ & $231.621$ & $28.20$ & $403.6$ & $10.7$ & $250.707$ & $28$ & $1$ & $1$ \\
		$400$ & $5$ & $7.90$ & $157.3$ & $5.6$ & $141.453$ & $30.80$ & $378.3$ & $10.3$ & $148.225$ & $26$ & $2$ & $2$ \\
		$400$ & $7$ & $8.67$ & $173.3$ & $6.6$ & $105.463$ & $24.40$ & $341.7$ & $9.3$ & $107.701$ & $23$ & $6$ & $1$ \\
		\bottomrule
	\end{tabular}
\end{table}

\subsubsection{Mean--variance model on market data}
\label{subsec:num-mv}
We consider the cardinality-constrained
mean--variance model
\begin{equation}\label{eq:mv-model}
	\begin{aligned}
		\min_{x\in\mathbb R^n}\quad & \frac12x^\top\Sigma x \\
		\text{s.t.}\quad & \mu^\top x\geq\rho, \ e^\top x=1, \  x\geq0, \  \|x\|_0\leq\kappa,
	\end{aligned}
\end{equation} 
where the variance of the portfolio is minimized subject to a lower
bound $\rho$ on its expected return. In contrast to the previous
problem, the set $C$ now combines an inequality and an equality
constraint, since we set
$G(x)=(\mu^\top x-\rho,\ e^\top x-1)^\top$,
$C=\mathbb R_+\times\{0\}$, and $\mathcal Q_C(s)=\diag(s_1^2,0)$.
We test five real stock data sets with daily return rates, which are
listed in \Cref{tab:num-mv-data} and include NASDAQ100, FTSE100,
and three S\&P500 data sets covering different market conditions.

\begin{table}[tbp]
	\centering
	\caption{The tested stock data sets.}
	\label{tab:num-mv-data}
	\footnotesize\setlength{\tabcolsep}{3.5pt}
	\begin{tabular}{llccc}
		\toprule
		Data set & Market conditions & Stocks & Total time period
		& Out-of-sample time period\\
		\midrule
		NASDAQ100 & -- & $95$ & 01/02/2020--09/12/2023 & 08/03/2022--09/12/2023 \\
		FTSE100 & -- & $97$ & 01/02/2020--09/12/2023 & 07/19/2022--09/12/2023 \\
		S\&P500-D & Declining market & $490$ & 03/12/2020--10/12/2022 & 01/03/2022--10/12/2022 \\
		S\&P500-C & COVID-19 & $478$ & 12/06/2017--12/31/2020 & 02/03/2020--12/31/2020 \\
		S\&P500-R & Rising market & $497$ & 09/17/2021--09/30/2024 & 11/01/2023--09/30/2024 \\
		\bottomrule
	\end{tabular}
\end{table}

The test problems are generated by a rolling-window procedure, in
which the first $70\%$ of the periods of each data set form the initial
estimation window, whereas the remaining $T_{\mathrm{out}}$ periods
form the out-of-sample period in \Cref{tab:num-mv-data}. For each
$t=1,\ldots,T_{\mathrm{out}}$, the estimation window is shifted by
$t-1$ periods, $\Sigma$ and $\mu$ are taken as the sample covariance
matrix and the sample mean of the returns in this window, and problem
\eqref{eq:mv-model} is solved with $\kappa=10$ and
$\rho=(e^\top\mu-\max_i\mu_i)/n$. The resulting portfolio is then
evaluated on the return of the following period.
For CD, we set $\tau=10$, $\eta=0.5$, $c=10$, $\alpha=0.01$, and
$\beta_0=10$.

\begin{figure}[tbp]
	\centering
	\includegraphics[width=\textwidth]{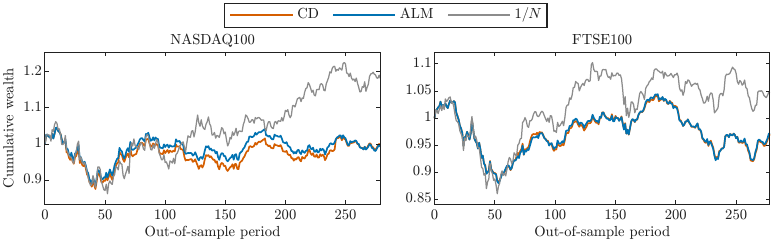}
	\caption{Cumulative wealth over the out-of-sample periods of
		NASDAQ100 and FTSE100 for the portfolios computed by CD and ALM.}
	\label{fig:num-mv}
\end{figure}

\begin{figure}[tbp]
	\centering
	\includegraphics[width=\textwidth]{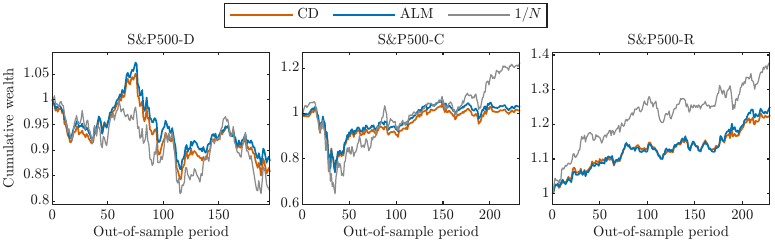}
	\caption{Cumulative wealth over the out-of-sample periods of the
		S\&P500 data sets for the portfolios computed by CD and ALM.}
	\label{fig:num-mv-regimes}
\end{figure}

Following \cite{wu2024sparse}, the computed portfolios are assessed out
of sample by the mean $Mean_p$ and the standard deviation $Std_p$ of
their realized returns, the Sharpe ratio $SR_p$, and the risk reduction
$RR_p$ with respect to the equally weighted portfolio, which is denoted
by $1/N$.
The corresponding results are summarized in
\Cref{tab:num-mv-out}, which reports, for each data set, the
averages of $Outer$, $Inner$, and $Time$ (in milliseconds) over all
windows, the number $N_{\mathrm{low}}$ of windows on which the
respective method attains the lower objective function value, and
these criteria. In addition, Figures~\ref{fig:num-mv}
and~\ref{fig:num-mv-regimes} show the cumulative wealth over the
out-of-sample periods.

\begin{table}[tbp]
	\centering
	\caption{Out-of-sample results of the rolling-window procedure for
		solving \eqref{eq:mv-model}.}
	\label{tab:num-mv-out}
	\footnotesize\setlength{\tabcolsep}{3.5pt}
	\begin{tabular}{lc l cccc cccc}
		\toprule
		Data set & $T_{\mathrm{out}}$ & Alg. & $Outer$ & $Inner$ & $Time$
		& $N_{\mathrm{low}}$ & $Mean_p$ & $Std_p$ & $SR_p$ & $RR_p$\\
		\midrule
		NASDAQ100 & $279$ & CD & $2.97$ & $499.4$ & $21.8$ & $155$ & $0.003$ & $0.0092$ & $-0.0125$ & $2.63$ \\
		 &  & ALM & $13.01$ & $1329.7$ & $28.2$ & $54$ & $0.002$ & $0.0091$ & $-0.0139$ & $2.59$ \\
		 &  & $1/N$ & -- & -- & -- & -- & $0.069$ & $0.0140$ & $0.0394$ & $1.00$ \\
		\addlinespace[2pt]
		FTSE100 & $278$ & CD & $3.28$ & $764.0$ & $27.0$ & $88$ & $-0.008$ & $0.0072$ & $-0.0315$ & $2.16$ \\
		 &  & ALM & $9.73$ & $440.5$ & $8.0$ & $63$ & $-0.008$ & $0.0071$ & $-0.0312$ & $2.15$ \\
		 &  & $1/N$ & -- & -- & -- & -- & $0.022$ & $0.0100$ & $0.0076$ & $1.00$ \\
		
		\addlinespace[2pt]
		S\&P500-D & $196$ & CD & $2.98$ & $464.1$ & $35.8$ & $134$ & $-0.073$ & $0.0104$ & $-0.0846$ & $2.32$ \\
		 &  & ALM & $15.82$ & $1622.9$ & $85.9$ & $49$ & $-0.061$ & $0.0101$ & $-0.0746$ & $2.34$ \\
		 &  & $1/N$ & -- & -- & -- & -- & $-0.091$ & $0.0147$ & $-0.0717$ & $1.00$ \\
		\addlinespace[2pt]
		S\&P500-C & $232$ & CD & $2.99$ & $544.5$ & $33.1$ & $30$ & $0.021$ & $0.0171$ & $0.0041$ & $3.19$ \\
		 &  & ALM & $6.61$ & $437.1$ & $20.6$ & $46$ & $0.029$ & $0.0173$ & $0.0083$ & $3.19$ \\
		 &  & $1/N$ & -- & -- & -- & -- & $0.118$ & $0.0249$ & $0.0418$ & $1.00$ \\
		\addlinespace[2pt]
		S\&P500-R & $229$ & CD & $3.05$ & $603.4$ & $49.8$ & $126$ & $0.091$ & $0.0057$ & $0.1352$ & $1.84$ \\
		 &  & ALM & $10.69$ & $1125.6$ & $61.7$ & $74$ & $0.099$ & $0.0059$ & $0.1441$ & $1.79$ \\
		 &  & $1/N$ & -- & -- & -- & -- & $0.143$ & $0.0077$ & $0.1679$ & $1.00$ \\
		
		\bottomrule
	\end{tabular}
\end{table}


Regarding the objective function values, CD attains the
lower value on $533$ of the $1214$ windows and ALM on $286$, while the
two methods tie on the remaining $395$ windows, $156$ of which belong
to S\&P500-C. In particular, CD attains the lower value more often on
four of the five data sets.
Regarding the numerical effort, CD requires about 3 outer
iterations on average on every data set, compared with between $6.61$
and $15.82$ for ALM, and, on NASDAQ100, S\&P500-D, and S\&P500-R, also
fewer inner iterations and less running time, the difference being most
pronounced on S\&P500-D, where the average running time is $35.8$~ms
for CD and $85.9$~ms for ALM.
These are the data sets on which the return constraint is active in
most windows. In fact, on the $776$ windows on which this
constraint is active (up to $10^{-4}$) at the point computed by CD,
CD requires $447.8$ inner iterations on average, compared with
$1309.8$ for ALM, whereas ALM requires fewer inner iterations on the
remaining $438$ windows.

Finally, since both methods address the same model, they yield similar
out-of-sample criteria. On the windows without a tie, the in-sample
variances of the portfolios computed by the two methods differ on
average by at most $1.2\%$, and the cumulative wealth curves in
Figures~\ref{fig:num-mv} and~\ref{fig:num-mv-regimes} remain close to
each other, with final wealth differing by at most $0.024$.
Hence, the savings in computational effort reported above are achieved
without a deterioration of the out-of-sample performance beyond these
small differences.
The sparse portfolios have a smaller
standard deviation than the equally weighted portfolio, with
$RR_p$ between $1.8$ and $3.2$, and, in the declining market of
S\&P500-D, they also lose less, with a final wealth of $0.86$ for CDiP and $0.88$ for ALM, compared
with $0.82$.

\subsection{MAXCUT problems}
\label{subsubsec:num-maxcut}
We consider an undirected graph with vertex set $V=\{1,\ldots,N\}$ and
denote $a_{ij}=a_{ji}$ as the weight of the edge between the vertices
$i$ and $j$. Since we allow zero weights, without loss of
generality, we can assume that the graph is complete.
Given a subset $S\subseteq V$ with complement $S^c$, the cut defined by
$S$ is the set of all edges with one end point in $S$ and the other one
in $S^c$, and its weight is $\sum_{i\in S,\,j\in S^c}a_{ij}$.
The MaxCut problem looks for a cut with maximum weight and is known to
be NP-hard.

In the following, let us denote $\mathbb S^N$ as the space of symmetric $N\times N$
matrices, equipped with the Frobenius inner product.
With $A:=(a_{ij})\in\mathbb S^N$ and the Laplacian matrix
$L:=\operatorname{Diag}(Ae)-A$, by \cite{jia2023augmented, goemans1995improved}, the MaxCut
problem can be reformulated as
\begin{equation}\label{eq:maxcut}
	\begin{aligned}
		\min_{W\in\mathbb S^N}\quad & f(W):=-\frac14\operatorname{tr}(LW),\\
		\text{s.t.}\quad & \diag(W)=e, \  W\succeq0, \ 
		\operatorname{rank}(W)\leq1,
	\end{aligned}
\end{equation}
since the feasible points of \eqref{eq:maxcut} are exactly the matrices
$W=ww^\top$ with $w\in\{-1,1\}^N$, and, for such a matrix,
$-f(W)=\frac14w^\top Lw$ is the weight of the cut defined by
$S=\{i\in V:w_i=1\}$.
Omitting the rank constraint yields the semidefinite relaxation of
\eqref{eq:maxcut}.
We take the set of positive semidefinite matrices of rank at most one as the complicated set, that
is, we set $G(W)=\diag(W)-e$, $C=\{0\}$, and
$D=\{W\in\mathbb S^N:W\succeq0,\ \operatorname{rank}(W)\leq1\}$.
A projection of $W\in\mathbb S^N$ onto $D$ is given by
$\max\{\lambda,0\}\,vv^\top$, where $\lambda$ is a largest eigenvalue of
$W$ and $v$ is a corresponding unit eigenvector, so that its
computation does not require a full spectral decomposition.
For the constraint dissolving mapping, we choose
$\mathcal Q_C\equiv0$ and $\mathcal Q_D(W)[Y]=c\,(WY+YW)$.

We use the rudy collection of the Biq Mac library  \href{http://biqmac.aau.at/biqmaclib.html}{http://biqmac.aau.at/\allowbreak biqmaclib.html} as the test database,
which consists of $130$ randomly generated graphs with $N\in\{60,80,100\}$
vertices, edge densities between $0.1$ and $0.99$, and integer edge
weights in $[-10,10]$ and for which the optimal cut values are known.
In order to obtain a good starting point, we take as $\widetilde W^0$ a
solution of the semidefinite relaxation of \eqref{eq:maxcut}, which we
compute by ALM with $D$ replaced by the cone of positive semidefinite
matrices.
For CD, we set $\tau=10$, $\eta=0.5$, $c=1$, $\alpha=0.3$, and
$\beta_0=5$. From the final iterate $W$ of each method, we obtain a cut from the
signs of the entries of an eigenvector of $W$ corresponding to its
largest eigenvalue.

\begin{figure}[tbp]
	\centering
	\includegraphics[width=0.87\textwidth]{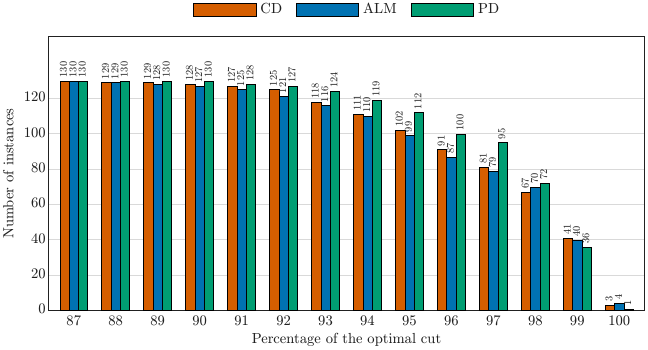}
	\caption{Number of the $130$ rudy instances on which CD, ALM, and PD
		attain a cut of at least $p\%$ of the optimal value.}
	\label{fig:num-maxcut}
\end{figure}

\Cref{fig:num-maxcut} summarizes
the cut quantity in a condensed way: for each instance, we take the
quotient $f_{\text{final}}/f_{\text{opt}}\in [0,1]$ as quantity measurement of the cut, where $f_{\text{final}}$ is the final function value obtained by each of the three methods and $f_{\text{opt}}$ is the optimal cut value.
\Cref{fig:num-maxcut} then reports, for $p=87,\ldots,100$, the number of
instances on which each method attains at least $p\%$ of the optimal
value.
For $p=88,\ldots,98$, PD attains at least $p\%$ of the optimal value on
more instances than CD and ALM, and its cuts attain on average
$97.50\%$ of the optimal values, compared with $97.12\%$ for CD and
$96.98\%$ for ALM. 

However, PD requires more outer and inner iterations and longer running times compared with CD and ALM. This is reflected in \Cref{tab:num-maxcut}, which reports
the average numbers of outer iterations ($Outer$) and of cumulative
inner iterations ($Inner$), the average running time in seconds
($Time$), and the average weight of the computed cuts in percent of the
optimal values ($Cut$).
\begin{table}
	\centering
	\caption{Results of CD, ALM, and PD for solving \eqref{eq:maxcut} with the rudy instances. Each row except the last summarizes the $10$
		instances of one group of the rudy collection.  The last row reports the average results over all instances listed above.}
	\label{tab:num-maxcut}
	\footnotesize\setlength{\tabcolsep}{2pt}
	\begin{tabular}{l cccc cccc cccc}
		\toprule
		& \multicolumn{4}{c}{CD} & \multicolumn{4}{c}{ALM} & \multicolumn{4}{c}{PD}\\
		\cmidrule(lr){2-5}\cmidrule(lr){6-9}\cmidrule(lr){10-13}
		Group
		& $Outer$ & $Inner$ & $Time$ & $Cut$
		& $Outer$ & $Inner$ & $Time$ & $Cut$
		& $Outer$ & $Inner$ & $Time$ & $Cut$\\
		\midrule
		\texttt{g05\_60} & $2.00$ & $27.6$ & $0.026$ & $99.00$ & $11.60$ & $227.4$ & $0.067$ & $99.44$ & $130.6$ & $103228$ & $30.7$ & $98.93$ \\
		\texttt{g05\_80} & $2.00$ & $34.9$ & $0.026$ & $99.07$ & $16.50$ & $259.7$ & $0.102$ & $98.87$ & $136.5$ & $149546$ & $61.0$ & $99.02$ \\
		\texttt{g05\_100} & $4.00$ & $42.4$ & $0.045$ & $99.07$ & $20.40$ & $365.9$ & $0.219$ & $98.98$ & $141.7$ & $202753$ & $122.8$ & $98.99$ \\
		\texttt{pm1s\_80} & $1.40$ & $25.3$ & $0.027$ & $93.37$ & $10.60$ & $220.7$ & $0.095$ & $95.63$ & $123.9$ & $112677$ & $45.3$ & $95.73$ \\
		\texttt{pm1d\_80} & $2.00$ & $61.2$ & $0.040$ & $96.31$ & $16.70$ & $320.3$ & $0.125$ & $92.77$ & $136.1$ & $157439$ & $63.8$ & $96.81$ \\
		\texttt{pm1s\_100} & $3.00$ & $61.0$ & $0.069$ & $95.90$ & $12.10$ & $288.6$ & $0.175$ & $94.41$ & $128.0$ & $150592$ & $89.7$ & $96.89$ \\
		\texttt{pm1d\_100} & $2.00$ & $55.8$ & $0.057$ & $96.00$ & $20.50$ & $376.5$ & $0.219$ & $95.87$ & $140.0$ & $227674$ & $135.8$ & $97.03$ \\
		\texttt{w01\_100} & $2.00$ & $77.1$ & $0.080$ & $96.29$ & $19.70$ & $382.0$ & $0.224$ & $97.65$ & $147.6$ & $289699$ & $173.0$ & $96.60$ \\
		\texttt{pw01\_100} & $4.00$ & $55.4$ & $0.058$ & $97.91$ & $14.10$ & $396.9$ & $0.233$ & $98.97$ & $151.9$ & $264581$ & $159.3$ & $97.95$ \\
		\texttt{w05\_100} & $4.00$ & $110.7$ & $0.113$ & $95.35$ & $17.20$ & $436.8$ & $0.259$ & $95.46$ & $155.7$ & $344831$ & $205.9$ & $95.68$ \\
		\texttt{pw05\_100} & $4.00$ & $51.3$ & $0.051$ & $99.09$ & $19.50$ & $378.5$ & $0.220$ & $99.12$ & $161.5$ & $279766$ & $170.8$ & $98.93$ \\
		\texttt{w09\_100} & $4.00$ & $124.7$ & $0.128$ & $95.87$ & $17.80$ & $408.5$ & $0.239$ & $94.31$ & $158.7$ & $403367$ & $239.2$ & $95.77$ \\
		\texttt{pw09\_100} & $6.00$ & $63.1$ & $0.063$ & $99.32$ & $23.20$ & $427.1$ & $0.248$ & $99.27$ & $164.8$ & $291427$ & $178.3$ & $99.19$ \\
		\midrule
		All & $3.11$ & $60.8$ & $0.060$ & $97.12$ & $16.92$ & $345.3$ & $0.187$ & $96.98$ & $144.4$ & $229045$ & $128.9$ & $97.50$ \\
		\bottomrule
	\end{tabular}
\end{table}
In fact, PD requires on average $128.9$~s running time per instance, which is about
$690$ times that of ALM and about $2100$ times that
of CD. On every instance, its running time exceeds those of CD and
ALM by more than two orders of magnitude.

\begin{figure}[tbp]
	\centering
	\includegraphics[width=\textwidth]{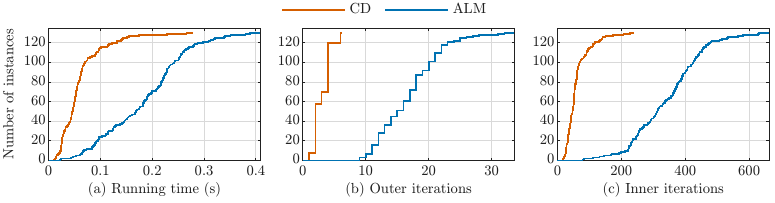}
	\caption{Number of the $130$ rudy instances that CD and ALM solve
		within a given running time~(a), number of outer iterations~(b),
		and number of cumulative inner iterations~(c).}
	\label{fig:num-maxcut-effort}
\end{figure}

Therefore, we merely compare the numerical performance of CD and ALM in the form of \cite{dolan2002benchmarking}. Specifically, \Cref{fig:num-maxcut-effort} shows the number of successfully-solved instances within a given running time, number of outer iterations, and number of cumulative
inner iterations. Overall, CD consistently outperforms ALM: it solves more instances faster and with fewer iterations.

\section{Conclusion}\label{Sec:conclusion}

In this manuscript, we developed a constraint dissolving method for solving optimization problems with constraints $G(x)\in C$ and $x\in D$, which, to the best of our knowledge, is the first constraint dissolving type method applicable to this general setting. More specifically, we introduced an auxiliary variable to decouple the constraint $G(x)\in C$ into an equality constraint and a set constraint, and constructed a constraint dissolving mapping acting jointly on the original and auxiliary variables. Under mild conditions, every accumulation point generated by the proposed algorithm is guaranteed to be an M-stationary point. Moreover, by a specific choice of the constraint dissolving mapping, every accumulation point is feasible. Numerical experiments on Rosenbrock, portfolio, and MAXCUT problems demonstrated the effectiveness of the proposed method in comparison with the augmented Lagrangian method and the penalty decomposition method, which, to the best of our knowledge, are the state-of-the-art algorithms for solving optimization problems with general constraints.

\vspace{3mm}
\noindent
{\bf Data Availability.} 
The five real-world stock datasets used in Section~5.2.2, including NASDAQ100, FTSE100, and three S\&P500 datasets, are available from Haonan Shen (15905823189@163.com) upon reasonable request.

\appendix

\section{Proof of \Cref{prop:block-cdm}}\label{app:proof-block-cdm}
\begin{proof}
By \Cref{ass:P} and \Cref{ass:projective}, $G$, $\nabla G$, $\mathcal Q_D$ and $\mathcal Q_C$ are continuously differentiable, and thus $\mathcal A(x,s)$ defined in \eqref{eq:cdm-def} is also continuously differentiable. Also, for all $(x,s)\in \mathcal F_s$, we have $H(x,s)=0$, and thus $\mathcal A(x,s)=(x,s)$. 

In view of \eqref{eq:cdm-def}, we have
\begin{equation}\label{Eq:der_A}
    \nabla \mathcal A(x,s) = I_{n+m} - \mathcal D J(x,s) H(x,s)-J(x,s) \nabla H(x,s), \  \forall \ (x,s)\in K,
\end{equation}
where $\mathcal D J(x,s)$ denotes the Jacobian of $J$ at $(x,s)$. In particular,
\begin{equation*}
    \nabla \mathcal A(x,s)= I_{n+m} - J(x,s) \nabla H(x,s), \ \forall \ (x,s)\in \mathcal F_s. 
\end{equation*}
Moreover, for any $(x,s)\in \mathcal F_s$, we have $M(x,s)=M_0(x,s)$. Then \eqref{eq:block-matrices} and \eqref{eq:correction-matrix} yield
\begin{equation*}\label{Eq:gradH+J} 
    \nabla H(x,s)J(x,s) = \nabla H(x,s) \mathcal Q(x,s) \nabla H(x,s)^ \top \left(\nabla H(x,s) \mathcal Q(x,s) \nabla H(x,s)^ \top \right) ^{-1} 
     = I_{m}. 
\end{equation*}
Consequently, one has
\begin{equation*}
    \nabla \mathcal A(x,s)^\top \nabla H(x,s)^\top= \nabla H(x,s)^\top - \nabla H(x,s)^\top J(x,s)^\top \nabla H(x,s)^\top =0,
\end{equation*}
which means \Cref{def:cdm} \ref{item:Ac} holds. It remains to verify \Cref{def:cdm} \ref{item:Ad} and \ref{item:Ae}.

Motivated by \eqref{Eq:der_A}, we set $\mathcal R_{\mathcal A}(x,s):=-\mathcal D J(x,s)H(x,s).$ 
By \Cref{ass:P} and \Cref{ass:projective}, we know that $\mathcal R_{\mathcal A}$ is continuous. For any $(x,  s)\in \mathcal F_s$, we have $H(x, s)=0$ and thus $\mathcal R_{\mathcal A} (x, s)=0$. It remains to verify \ref{item:Ae}.

For any $(x,  s)\in \mathcal F_s$, choose a constant $\omega_{x, s}>0$ arbitrarily. For any $(y,t)\in K$ satisfying $\|(y,t)-(x, s)\| \leq \omega_{x, s}$, by \Cref{ass:projective} \ref{item: Q}, one has $\mathcal Q(y,t) \mathcal N_K^{\text{lim}}(y,t)=\{0\}$ and thus
\begin{equation*}
\begin{aligned}
    &\left( \nabla \mathcal A(y,t)- I_{n+m} - \mathcal R_{\mathcal A}(y,t)\right)^ \top \mathcal N_K^{\text{lim}}(y,t)= \left( -J(y,t)\nabla H(y,t) \right)^\top \mathcal N_K^{\text{lim}}(y,t)\\
    &=-\nabla H(y,t)^\top M(y,t)^{-1} \nabla H(y,t) \mathcal Q(y,t) \mathcal N_K^{\text{lim}}(y,t) =\{0\}.
\end{aligned}  
\end{equation*}
This completes the proof.
\end{proof}

\bibliographystyle{tfs}
\bibliography{bibliography-3}

\end{document}